\documentclass[12pt,twoside,leqno]{article}
\usepackage{amsmath}
\usepackage{amssymb}
\usepackage{amsxtra}
\usepackage{amscd}
\usepackage{amsthm}
\usepackage[mathscr]{eucal}

\theoremstyle{plain}
\newtheorem{thm}[subsection]{Theorem}
\newtheorem{thmstar}[subsection]{Theorem*}
\newtheorem{prop}[subsection]{Proposition}
\newtheorem{propstar}[subsection]{Proposition*}
\newtheorem{cor}[subsection]{Corollary}
\newtheorem{lem}[subsection]{Lemma}

\newtheorem{sblem}[subsubsection]{Lemma}

\theoremstyle{definition}
\newtheorem{defn}[subsection]{Definition}

\newtheorem{rem}[subsection]{Remark}
\newtheorem{para}[subsection]{}

\newtheorem{sbpara}[subsubsection]{}

\newenvironment{pf}{\proof[\proofname]}{\endproof}

\begin{document}
\title
{Logarithmic structures of Fontaine-Illusie.\  II.\  \\
--- Logarithmic flat topology.
}
\author{
Kazuya Kato
}
\date{}
\maketitle

\newcommand\Cal{\mathcal}
\newcommand\define{\newcommand}
\renewcommand\bold{\Bbb}

\define\bG{\bold G}
\define\bZ{\bold Z}
\define\bC{\bold C}
\define\bR{\bold R}
\define\bQ{\bold Q}
\define\bN{\bold N}
\define\bP{\bold P}
\define\cl{\mathrm{cl}}%
\define\fil{\mathrm{fil}}%
\define\fl{\mathrm{fl}}%
\define\gp{\mathrm{gp}}%
\define\fs{\mathrm{fs}}%
\define\an{\mathrm{an}}%
\define\et{\mathrm{\acute{e}t}}%
\define\mult{\mathrm{mult}}%
\define\sat{\mathrm{sat}}%
\renewcommand\int{\mathrm{int}}%
\define\Ker{\mathrm{Ker}\,}%
\define\Coker{\mathrm{Coker}\,}%
\define\Hom{\operatorname{\mathrm{Hom}}}%
\define\Aut{\operatorname{\mathrm{Aut}}}%
\define\Mor{\operatorname{\mathrm{Mor}}}%
\define\rank{\mathrm{rank}\,}%
\define\gr{\mathrm{gr}}%
\define{\cHom}{\operatorname{\mathcal{H}\mathit{om}}}
\define{\HOM}{\cHom}
\define{\cExt}{\operatorname{\mathcal{E}\mathit{xt}}}
\define{\cMor}{\operatorname{\mathcal{M}\mathit{or}}}
\define\cO{\Cal O}
\renewcommand{\O}{\cO}
\define\cS{\Cal S}
\define\cM{\Cal M}
\define\cG{\Cal G}
\define\cH{\Cal H}
\define\cE{\Cal E}
\define\cF{\Cal F}

\newcommand{\N}{{\mathbb{N}}}
\newcommand{\Q}{{\mathbb{Q}}}
\newcommand{\Z}{{\mathbb{Z}}}
\newcommand{\R}{{\mathbb{R}}}
\newcommand{\C}{{\mathbb{C}}}
\newcommand{\ol}[1]{\overline{#1}}
\newcommand{\respect}{\rightsquigarrow}
\newcommand{\compatible}{\leftrightsquigarrow}
\newcommand{\upc}[1]{\overset {\lower 0.3ex \hbox{${\;}_{\circ}$}}{#1}}
\newcommand{\Gmlog}{\bG_{m, \log}}%
\newcommand{\Gm}{\bG_m}%
\newcommand{\Ga}{\bG_a}%
\newcommand{\ep}{\varepsilon}
\newcommand{\pe}{\frak p}
\newcommand{\Spec}{\operatorname{Spec}}
\newcommand{\val}{{\mathrm{val}}}
\newcommand{\lcf}{{\mathrm{lcf}}}
\newcommand{\Et}{{\mathrm{\acute{E}t}}}
\renewcommand{\b}{\langle\,,\,\rangle}%
\newcommand{\n}{\operatorname{naive}}
\newcommand{\bs}{\operatorname{\backslash}}
\newcommand{\Lie}{\operatorname{Lie}}

\renewcommand{\emph}{\it}
\newcommand{\Map}{\operatorname{Map}}

\begin{abstract} This is Part II of the author's paper Logarithmic structures of Fontaine-Illusie. We discuss log flat topology and log flat descent. We study the first log flat cohomology
$H^1(X_{\log, \fl}, G)$ for various sheaves of groups $G$, for example, $G=GL_n$, finite flat commutative group schemes, the log multiplicative group $M^{\gp}$, etc. 
\end{abstract}

{\bf Contents}

\medskip

\S\ref{sec:basic}. Basic facts.  

\S\ref{sec:flat}. Logarithmic flat topology. 

\S\ref{sec:des1}. Descent theory.\ I.

\S\ref{sec:coh}. Cohomology. 

\S\ref{sec:Hilb}. Hilbert 90.

\S\ref{sec:free}. Locally free modules. 

\S\ref{sec:des2}. Descent theory.\ II.\ Logarithmic flat descent

\S\ref{sec:des3}. Descent theory.\ III.

\S\ref{sec:torsor}.  Torsors.

\S\ref{sec:pi1}.  Logarithmic fundamental groups.

\bigskip

\section*{Introduction.}
  This is a continuation of the paper \cite{log} on the foundation of 
log geometry in the sense of Fontaine-Illusie. 
  Here we discuss mainly log flat topologies, especially log flat 
descent theory.  

  This paper was started 
around 1991, and was circulated 
as an incomplete preprint for a long time.  
  Since then, some contents of this paper have been reproduced by several 
authors with proofs 
(\cite{H}, \cite{KS}, \cite{Ni}, \cite{N:qs}, \cite{O}, ...). 
  A.\ Moriwaki \cite{M} also studied flat descents in the category of log schemes. 
  In the parts which were incomplete in the circulated preprint, 
we sometimes referred to these papers instead of completing the original proofs. 
  In particular, the author does not claim the results with $*$ 
(i.e., two theorems \ref{thm:des2}, \ref{thm:ketdes} and one proposition \ref{prop:zar}) 
are his results.
  Since the paper is already referred to in many published works, in the other parts, we preferred to preserve the original, circulated form.
  In both parts, 
we tried to preserve the original numberings of definitions and propositions.

  The author wishes to express his special thanks to Chikara Nakayama who helped him a lot in the completion of this paper. He is also 
 thankful to Luc Illusie, Takeshi Saito, 
Takeshi Kajiwara, and Takeshi Tsuji for helpful discussions. 

The author is partially supported by NSF Award 1601861. 

\section{Basic facts.}\label{sec:basic}

\bigskip

  In this section, we review basic facts about logarithmic 
schemes in \cite{log}, and introduce ``fs log schemes" and 
``log flat morphisms." 

\begin{para}
  In this paper, a monoid is always assumed to be commutative 
and has a unit element. 
  A homomorphism of monoids is assumed to respect the unit 
elements. 
  The semi-group law of a monoid is usually written 
multiplicatively. 

  For a monoid $P$, we denote by $P^{\times}$ the group of all 
invertible elements of $P$, and by $P^{\gp}$ the group hull 
$P^{-1}P =\{a^{-1}b \ ; \ a, b \in P\}$ of 
$P$. 
\end{para}

\begin{para}
{\emph Pre-logarithmic structures, logarithmic structures, 
logarithmic schemes}. 
  A pre-log structure on a scheme $X$ is a sheaf of monoids 
$M$ on the \'etale site $X_{\et}$ endowed with a homomorphism 
of sheaves of monoids 
$$ \alpha \colon M \to \cO_X,$$ 
where $\cO_X$ is regarded as a monoid for the multiplicative 
law. 

  A pre-log structure $M$ is called a log structure if 
$$\alpha^{-1}(\cO_X^{\times}) \overset \cong \to \cO_X^{\times} 
\qquad \mathrm{via} \ \alpha.$$ 

  A log scheme is a scheme endowed with a log structure. 
  A morphism of log schemes is defined in a natural way. 
  For a log scheme $X$, the log structure of $X$ is usually denoted 
by $M_X$. 
\end{para}

\begin{para}
  {\emph The log structure associated to a pre-log structure}. 
  For a scheme $X$ and a pre-log structure $M \overset \alpha 
\to \cO_X$, the log structure $M^{\sim}$ associated to $M$ is 
defined to be the pushout of the diagram 
$$M \leftarrow \alpha^{-1}(\cO_X^{\times}) \overset \alpha 
\to \cO_X^{\times}$$ 
in the category of sheaves of monoids, which is endowed with 
the homomorphism $M^{\sim} \to \cO_X$ induced by $\alpha$ and 
by $\cO_X^{\times} \overset \subset \to \cO_X$. 
\end{para}

\begin{para}
  {\emph The inverse image of a log structure}. 
  Let $f\colon X \to Y$ be a morphism of schemes, and let $M$ 
be a log structure on $Y$. 
  Then the inverse image $f^*M$ of $M$ is defined to be the 
log structure on $X$ associated to the pre-log structure 
$f^{-1}(M)$ endowed with $f^{-1}(M) \to f^{-1}(\cO_Y) \to \cO_X$. 
  Here $f^{-1}(M)$ denotes the sheaf-theoretic inverse image of 
$M$. 
\end{para}

\begin{defn}
\label{defn:fs}
  (1)  We say a monoid $P$ is {\emph integral} if the canonical 
map $P \to P^{\gp}$ is injective (that is, if $ab=ac$ $(a, b, c 
\in P)$ implies $b=c$). 

  We regard an integral monoid $P$ as a submonoid of $P^{\gp}$. 

  (2)  We say a monoid $P$ is {\emph saturated} if $P$ is integral 
and satisfies 
the following condition: 
  If $a \in P^{\gp}$ and $a^n\in P$ for some $n\ge1$, then 
$a \in P$. 

  (3)  We call a finitely generated integral monoid a {\emph fine 
monoid}.  
  We call a finitely generated saturated monoid an {\emph fs 
monoid}. 
\end{defn}

\begin{para}
{\emph Fine log schemes and fs log schemes}. 
  We call a log scheme $X$ a fine (resp.\ 
an fs) log scheme if the following condition 
is satisfied: 
  \'Etale locally on $X$, there exists a fine (resp.\ an fs) 
monoid $P$ and a homomorphism $\alpha\colon P \to \cO_X$ such 
that $M_X$ is isomorphic to the log structure associated to the 
constant sheaf $P$ on $X$ regarded as a pre-log structure with 
respect to $\alpha$. 

  An fs log scheme is fine. 
  A fine log scheme $X$ is fs if and only if $M_X$ is a sheaf 
of saturated monoids. 

  For a fine (resp.\ fs) log scheme $X$, the stalk 
$(M_X/\cO_X^{\times})_{\overline x}$ is a fine monoid (resp.\ 
an fs monoid) for any $x\in X$. 
\end{para}

\begin{para}
\label{chart}
{\emph Charts}. 
  (1) A chart of a fine (resp.\ an fs) log scheme $X$ is a pair 
$(P, h)$, where $P$ is a fine (resp.\ an fs) monoid and $h$ is 
a homomorphism $P \to M_X$ satisfying the following 
condition:  
  Let $P^{\sim}$ be the log structure associated to $P \to \cO_X$ 
induced by $h$. 
  Then the induced map $P^{\sim}\to M_X$ is an isomorphism. 

  A chart for $X$ exists \'etale locally on $X$. 

  (2) Let $f\colon X \to Y$ be a morphism of fine (resp.\ fs) 
log schemes. 
  A chart of $f$ is a diagram $M_Y\overset \alpha \leftarrow 
P \to Q \overset \beta \to M_X$ such that $\alpha$ is a chart 
for $Y$, $\beta$ is a chart for $X$, and the diagram 
$$\begin{CD}
  P @>>> Q \\
 @VVV @VVV \\ 
  f^{-1}(M_Y) @>>> M_X
\end{CD}
$$
is commutative. 

  A chart of $f$ exists \'etale locally on $X$ and $Y$. 
  This is proved in \cite{log} in the fine case and the fs case is proved similarly.

  For a morphism $f\colon X \to Y$ and a fine (resp.\ an fs) log 
structure $M$ on $Y$, the inverse image $f^*M$ is also a fine 
(resp.\ an fs) log structure.  A chart of $M$ gives naturally 
a chart of $f^*M$. 
\end{para}

\begin{para}
{\emph Finite inverse limits}. 
  The category of fine (resp.\ fs) log schemes has finite inverse 
limits. 
  For a finite inverse system $\Sigma$, its 
inverse limit is a fine (resp.\ an fs) log scheme whose underlying 
scheme is finite and of finite presentation 
over the inverse limit of the underlying diagram 
of schemes of $\Sigma$. 
  To see this, it is enough to show it for fiber products. 
  The fiber product of $\Sigma\colon X_1 \overset {f_1} \to X_0 \overset {f_2} \leftarrow X_2$ is constructed \'etale locally on $X_i$'s as follows. 
  \'Etale locally, take charts $\alpha_i \colon P_i \to M_{X_i}$ of 
$X_i$ and charts $M_{X_0} \overset {\alpha_0} \leftarrow 
P_0 \overset {h_j} \to P_j \overset 
{\alpha_j} \to M_{X_j}$ 
of $f_j$ ($j=1, 2$) 
(this is possible by \ref{chart}).
  Let $X$ be the fiber product of the underlying diagram of 
schemes of $\Sigma$ and let $P$ be the pushout of the 
diagram $P_1 \overset {h_1} \leftarrow P_0 \overset {h_0} \to P_2$ in the category of 
monoids. 
  Then we have a canonical homomorphism $P \to \cO_X$. 
  Let $P^{\int}$ be the image of $P \to P^{\gp}$ (resp.\ 
$P^{\sat} = \{x \in P^{\gp}\ ; \ x^n \in \mathrm{Image}\, 
(P \to P^{\gp})$ for some $n \ge1\})$. 
  Then the fiber product 
of $\Sigma$ in the category of fine 
(resp.\ fs) log schemes is the scheme $T = X \otimes_{\Z[P]}
\Z[P^{\int}]$ (resp.\  $X \otimes_{\Z[P]}
\Z[P^{\sat}]$) endowed with the log structure associated to 
$P^{\int}$ (resp.\ $P^{\sat}$) $\to \cO_T$. 
\end{para}

\begin{para}
{\emph Advantages of fs log schemes}. 
  In \S\ref{sec:flat}--\S\ref{sec:pi1} of this paper, we 
consider fs log schemes but do not consider fine log schemes, 
though we considered fine log schemes in \cite{log}. 
  The category of fs log schemes has, for example, the following 
advantages. 

\begin{sbpara}
  The category of fs log schemes has much more group objects 
than the category of fine log schemes. 
  This point is explained in the preprint \cite{K2}. 
\end{sbpara}

\begin{sbpara}
  The exactness of the sequence 
$$0\to\Gamma(X,\bZ/n(1)) \to \Gamma(X,M_X^{\gp}) \overset n \to 
\Gamma(X, M_X^{\gp})$$
$(n\in \bZ, n\not=0)$, which plays important roles in log algebraic 
geometry (see \S\ref{sec:coh}), holds for an fs log scheme $X$, but 
does not hold in general for a fine log scheme $X$. 
\end{sbpara}

\begin{sbpara}
  ``Log Galois theory" works well for fs log schemes (see 
\S\ref{sec:pi1}) but not for fine log schemes. 
  For example, let $A$ be a discrete valuation ring, $n\ge2$ an 
integer which is invertible in $A$, $\pi$ a prime element of $A$, 
$B=A[\pi^{1/n}]$, and assume that $A$ contains a primitive $n$-th 
root of $1$.
  Let $G=\Aut_A(B)\cong \bZ/n$. 
  Endow $Y=\Spec(A)$ (resp.\ $X=\Spec(B)$) with the standard log 
structure which is associated to $\bN \to A; 1 \mapsto \pi$ 
(resp.\ $\bN \to B; 1\mapsto \pi^{1/n})$. 
  Then, 
$$G\times X \overset \cong \to X \times_YX; (g, x)\mapsto (x, gx)$$
(an isomorphism as in the \'etale Galois theory in classical 
algebraic geometry) holds in the category of fs log schemes, but 
does not hold in the category of fine log schemes. 
  (Indeed, the underlying scheme of the fiber product $X\times_YX$ 
in the category of fine log schemes is the same as the fiber 
product of the underlying schemes, and hence is connected (not the 
disjoint union of $\sharp(G)$-copies of $X$).) 
\end{sbpara}
\end{para}

\begin{para}
\label{defn:lf}
{\emph Log flatness, log smoothness, log \'etaleness}.
  Let $f\colon X\to Y$ be a morphism in either the category 
  of fine log schemes or the category of 
  fs log 
schemes. 
  We say $f$ is log flat (resp.\ log smooth, resp.\ log \'etale) 
if classical fppf (classical \'etale, classical \'etale) locally 
on $X$ and on $Y$, there exists a chart $M_Y \leftarrow P \to 
Q \to M_X$ of $f$ satisfying the following two conditions 
(i) and (ii). 

(i)  The induced map $P^{\gp}\to Q^{\gp}$ is injective (resp.\ 
is injective and the order of the torsion part of its cokernel is 
invertible on $X$, resp.\ is injective and its cokernel is finite 
with an order invertible on $X$). 

(ii) The induced morphism of schemes $X \to Y\times_{\bZ[P]}\bZ[Q]$ 
is flat (resp.\ smooth, resp.\ \'etale) in the classical sense. 

  A morphism of fs log schemes is log flat (resp.\ log smooth, 
resp.\ log \'etale) if and only if it is so as a morphism of fine 
log schemes.

  Log flat (resp.\ log smooth, resp.\ log \'etale) morphisms are 
stable under compositions and base changes.
  The stability for base changes is seen easily, and the stability 
of log smoothness and that of log \'etaleness for compositions are 
seen from the following infinitesimal characterization of log 
smoothness and that of log \'etaleness, respectively.   
 
   A morphism $f\colon X \to Y$ of fine log schemes is log smooth 
(resp.\ log \'etale) if and only if the following condition is 
satisfied: 

\noindent 
$(*)$  For any diagram of fine log schemes including $f$ 
$$
\begin{CD}
  T @>t>> X \\
 @ViVV @VVfV \\ 
  S @>s>> Y
\end{CD}
$$
such that $i^*(M_S) \overset \cong \to M_T$ and such that the 
underlying scheme of $T$ is (via $i$) a closed subscheme of the 
underlying scheme of $S$ defined by a nilpotent quasi-coherent 
ideal of $\cO_S$, then there exists classical \'etale locally on $S$ 
a morphism (resp.\ there exists a unique morphism) $g\colon S 
\to X$ such that $g\circ i = t$ and $f\circ g = s$. 

  See \cite{log} (3.5) for the proof.

  A morphism $f\colon X \to Y$ of fs log schemes is log smooth 
(resp.\ log \'etale) if and only if the condition \newline 
$(*){}'$ the same as $(*)$ except that we assume that $T$ and $S$ are 
fs log schemes, 
\newline 
is satisfied. 

  This is proved in the same way as the case of fine log schemes. 
\end{para}

 We prove the stability of log flatness for compositions. 
  This stability is also proved in \cite{O} Corollary 4.12 (ii) by another method. 
  See also Ogus' book \cite{Og} Chapter IV, Proposition 4.1.2 (4).

\begin{lem}\label{chartfl} Let $f: X\to Y$ be a log flat morphism of fine (resp.  fs)  log schemes and let $\beta: P \to M_Y$ be a chart.  Then fppf  locally on $X$ and on $Y$ in the classical sense, there exists a chart $(P\overset{\beta}\to M_Y, Q\to M_X, P\to Q)$ of $f$ including $\beta$ satisfying the conditions (i) and (ii)  in the definition of log flatness in \ref{defn:lf}. 
\end{lem}

\begin{pf} Let $(M_Y\leftarrow P' \to  Q'\to M_X)$ be a chart of $f$ satisfying these conditions (i) and (ii) in the definition of log flatness. 

{\bf Claim 1.} We may assume that there is a homomorphism $P\to P'$ such that  $\beta: P\to M_Y$ factors as $P\to P'\to M_Y$ and such that the homomorphism $P^{\gp}\to (P')^{\gp}/(P')^\times$ is surjective. 

{\it Proof} of Claim 1.  
Fix $x\in X$, $y=f(x)\in Y$. By replacing $P'$ with the inverse image $P''$ of $M_{Y, \overline{y}}$ under 
$$P^{\gp} \times (P')^{\gp} \to M_{Y, \overline{y}}^{\gp}\;;\; (a,b)\mapsto ab,$$
and by replacing $Q'$ with  the pushout of 
$P''\leftarrow P'\to Q'$ in the category of fine (resp. fs) monoids, we obtain the situation stated in Claim 1.

Assume we have $P\to P'$ as in Claim 1. 

{\bf Claim 2.} We may assume that there are a finitely generated abelian group $H$ which contains $P^{\gp}$ as a subgroup and a homomorphism $H\to (Q')^{\gp}$ such that the  diagram 
$$\begin{matrix}  P^{\gp}&\to& H\\
\downarrow && \downarrow\\
(P')^{\gp}&\to& (Q')^{\gp}\end{matrix}$$
is commutative and is co-cartesian in the category of abelian groups.

{\it Proof} of Claim 2. Let $R=(Q')^{\gp}/(P')^{\gp}$ and let $S$ be the inverse image of $(P')^\times$ in $P^{\gp}$. Then $P^{\gp}/S\overset{\cong}\to (P')^{\gp}/(P')^\times$. Consider the commutative diagram of exact sequences
$$\begin{matrix}  \text{Ext}^1(R, S) &\to  &\text{Ext}^1(R, P^{\gp}) &\to& \text{Ext}^1(R, P^{\gp}/S)&\to& 0\\
\downarrow && \downarrow&& \downarrow\\
 \text{Ext}^1(R, (P')^\times) &\to & \text{Ext}^1(R, (P')^{\gp}) &\to& \text{Ext}^1(R, (P')^{\gp}/(P')^\times)&\to& 0. \end{matrix}$$
(Note that $\text{Ext}^2=0$ for $\Z$-modules.)  Let $b\in  \text{Ext}^1(R, (P')^{\gp})$  be the class of $(Q')^{\gp}$. Since the right vertical arrow is an isomorphism, this  diagram shows that there is $a\in  \text{Ext}^1(R, P^{\gp})$   such that $a'-b$, where $a'$ denotes the image of $a$ in  $\text{Ext}^1(R, (P')^{\gp})$, comes from an element $c$ of  $\text{Ext}^1(R, (P')^\times)$. Since $c$ dies in $\text{Ext}^1(R, I)$ for a divisible abelian group $I$ containing $(P')^\times$ as a subgroup, $c$ dies in  $\text{Ext}^1(R, G)$ for some finitely generated   abelian group $G$ containing $(P')^\times$ as a subgroup. Replace $P'$ by the pushout of $P'\leftarrow (P')^\times\to G$, $Q'$ by the pushout of $Q'\leftarrow (P')^\times\to G$, $X$ by $\Z[G]\otimes_{\Z[(P')^{\times}]} X$, and $Y$ by  $\Z[G]\otimes_{\Z[(P')^{\times}]} Y$. Then $a'=b$ in $\text{Ext}^1(R, (P')^{\gp})$. We have the desired  abelian group $H$ with an exact sequence $0\to P^{\gp}\to H \to R \to 0$ corresponding to $a\in \text{Ext}^1(R, P^{\gp})$.  
This proves Claim 2. 

Now let the situation be as in Claim 2, and let $Q$ be the inverse image of $Q'$ in $H$. Since $Q/S\overset{\cong}\to Q'/(P')^\times$, $Q\to M_X$ is a chart of $X$. We have the chart $M_Y \overset{\beta}\leftarrow P \to Q \to M_X$ of $f$ having the desired property. 
\end{pf}

Now we complete the proof of the stability of the log flatness for composition. Let $f: X\to Y$ and $g: Y\to Z$ be log flat morphisms of fine (resp. fs) log schemes. Let $M_Z\leftarrow P \to Q \overset{\beta}\to M_Y$ be a chart of $g$ which satisfies the conditions (i) and (ii) in the definition of log flatness. By \ref{chartfl}, we may assume that there is a chart $M_Y \overset{\beta}\leftarrow Q \to R \to M_X$ of $f$ including $\beta$ which satisfies the conditions (i) and (ii) of log flatness. Then the chart $M_Z \leftarrow P \to R\to M_X$ of $g \circ f$ satisfies the conditions (i) and (ii) of log flatness. 

\medskip

  A morphism $f\colon X \to Y$ of fine log schemes 
such that the homomorphism $f^*M_X \to M_Y$ is an isomorphism 
is log flat 
(resp.\ log smooth, resp.\ log \'etale) if and only if the 
underlying morphism of schemes is flat (resp.\ smooth, resp.\ \'etale).

\section{The logarithmic flat topology and the logarithmic \'etale topology.}\label{sec:flat}

\bigskip

\begin{defn}\label{defn:Kummer} 
  A homomorphism of fs monoids $h\colon P \to Q$ is said to be 
{\emph of Kummer type} if the following condition is satisfied: 
It is injective, and for any $a \in Q$, there exists $n\ge1$ such 
that $a^n \in h(P)$. 
\end{defn}

\begin{defn}
  We say a morphism $f\colon X \to Y$ of fs log schemes is {\emph of 
Kummer type} if for any $x \in X$, the homomorphism of fs monoids 
$M_{Y, \overline{y}}/\cO^{\times}_{Y, \overline{y}} \to 
M_{X, \overline{x}}/\cO^{\times}_{X, \overline{x}}$ with $y = f(x)$ is of 
Kummer type in the sense of \ref{defn:Kummer}. 
\end{defn} 

  A standard example is the following. 
  Let $Y$ be an fs log scheme with a chart $P \to M_Y$, and 
assume that we are given an fs monoid $Q$ and a homomorphism 
$P \to Q$ of Kummer type. 
  Let $X$ be any scheme over $Y\otimes_{\bZ[P]}\bZ[Q]$ 
endowed with the log structure associated to $Q \to \cO_X$. 
  Then $f\colon X \to Y$ is of Kummer type. 
 
  If $X= Y\otimes_{\bZ[P]}\bZ[Q]$ here, $f$ is log flat 
and surjective, and $X$ is a standard example of a covering of $Y$ 
for the log flat topology introduced in the following 
\ref{defn:topology}.
 
  In the log flat case of 
\ref{defn:lf}, 
if $f$ is of Kummer type and locally of finite presentation, we can take the chart satisfying the additional 
conditions that $P \to Q$ is of Kummer type 
and that $X \to Y\otimes_{\bZ[P]}\bZ[Q]$ is surjective. 
See \cite{INT} Proposition 1.3.

  Now we define the logarithmic flat topology and the logarithmic 
\'etale topology on an fs log scheme. 
  (In the sequel we neglect questions of universes, which can be 
treated as in the non-log case.) 

\begin{defn}\label{defn:topology} 
  For an fs log scheme $T$, we say a family of morphisms 
$\{f_i\colon U_i \to T\}_i$ of fs log schemes is a covering of 
$T$ for the log flat (resp.\ log \'etale) topology, if the following 
conditions (i) (resp.\ (i)${}'$) and (ii) are satisfied. 

(i) $f_i$ are log flat and of Kummer type, and the underlying 
morphisms of schemes of $f_i$ are locally of finite presentation.

(i)${}'$ $f_i$ are log \'etale and of Kummer type. 

(ii) $T=\underset i \bigcup f_i(U_i)$ (set theoretically). 

  Let $X$ be an fs log scheme and let (fs/$X$) be the category of fs 
log schemes over $X$. 
  We define the Grothendieck topology called 
the log flat (resp.\ log \'etale) topology on 
(fs/$X$) by taking coverings as above.
  We denote the site $(\fs/X)$ endowed with this topology by 
$X^{\log}_{\fl}$ (resp.\ $X^{\log}_{\et}$). 

  To see that these really give definitions of Grothendieck topologies, 
we need the following lemma whose (2)  is a case of \cite{NC:loget} (2.2.2). 
  (In there, it is proved that the assumption that $f$ is of Kummer 
type in \ref{lem:4thpoint} (2) can be replaced by the weaker assumption that $f$ is exact.)
\end{defn}

\begin{lem}\label{lem:4thpoint}
  Let $f\colon X\to Y$ be a morphism of fs log schemes of Kummer type, let $Y'$ be an fs log scheme 
  over $Y$, and let $f':X':=X\times_Y Y'\to Y'$. 
  
 $(1)$ The morphism $f'$ is of Kummer type. 
  
 $(2)$ Assume 
that 
$f$ is surjective. Then $f'$ is surjective.
  More strongly, for any $x \in X$ and $y\in Y'$ having the same 
image in $Y$, there exists $z \in X'$ with image $x$ in 
$X$ and $y$ in $Y'$. 
\end{lem}

  Here and in the rest of this paper, when we are discussing fs log 
schemes, the notation of the fiber product stands for the fiber 
product in the category of fs log schemes unless the contrary is 
explicitly stated.

\begin{pf}
  We may assume that $X, Y$ and $Y'$ have the same underlying scheme 
which is the spectrum of an algebraically closed field $k$. Take a submonoid $P'$ of $M_{Y'}$ such that
$M_{Y'}=k^\times \times P'$. Let $P$ be the inverse image of $P'$ in $M_Y$ under the homomorphism $M_Y\to M_{Y'}$. 
Then $M_Y=k^\times \times P$. We have a chart $M_Y \leftarrow P \to P' \to M_{Y'}$ of $Y'\to Y$.

Take a submonoid $Q$ of $M_X$ such that $M_X=k^\times \times Q$. Then the map $M_Y\to M_X$ is written as $k^\times \times P\to k^\times \times Q\;;\;(s,t) \mapsto (sh(t), b(t))$ for some homomorphism $h: P \to k^\times$, where $b$ is the composition $P\cong M_Y/k^\times \to M_X/k^\times \cong Q$. Since $k^\times$ is a divisible abelian group and $b$ is injective, $h$ extends to a homomorphism $Q\to k^\times$ which we still denote by $h$. We have a commutative diagram 
$$\begin{matrix}  P&\overset{b}\to & Q \\
\downarrow &&\downarrow \\
M_Y&\to &M_X\end{matrix}$$
where the left vertical arrow is the inclusion map  and the right vertical arrow is $y \mapsto h(y)y$ which we denote by $\gamma$. We have a chart $M_Y \leftarrow P \overset{b}\to Q \overset{\gamma}\to M_X$ of $X\to Y$. 

Let $Q'$ be the pushout of $P'\leftarrow P \to Q$ in the category of fs monoids. Then the fiber product $X'$ is $\Spec(R)$ where $R=k[Q'] \otimes_{k[P']} k$ with the log structure associated to $Q'\to R$. Let $\Delta$ be the torsion subgroup of $(Q')^{\gp}$. Then $\Delta\subset Q'$ and the map $P'\to Q'/\Delta$ is injective. For any $a\in Q'\smallsetminus \Delta$, there is some $n\geq 1$ such that $a^n$ comes from $P'\smallsetminus \{1\}$. Hence there is a nilpotent ideal $I$ of $R$ and a nilpotent ideal $J$ of $k[\Delta]$ such that the map $k[\Delta]\to R$ induces $k[\Delta]/J\overset{\cong}\to R/I$. Hence $X'$ is not empty and $Q'/\Delta\overset{\cong}\to (M_{X'}/\cO^\times_{X'})_x$ for any $x\in X'$. This proves (1) and (2).  
\end{pf}

  We prove the following three propositions \ref{prop:open}, \ref{prop:exact}, \ref{prop:kummer} concerning coverings 
for the logarithmic flat topology. 
  First, the following \ref{prop:open} says that the family 
$\{f_i(U_i)\}_i$ in \ref{defn:topology} (ii) is an open covering 
of $T$. 

\begin{prop}\label{prop:open}
  Let $f\colon X \to Y$ be a morphism of fs log schemes. 
  Assume that 
$f$ is log flat and of Kummer type, and the underlying morphism 
of schemes of $f$ is locally of finite presentation. 
  Then $f$ is an open map. 
\end{prop}

\begin{prop}\label{prop:exact}
  Let $f\colon X \to Y$ be a morphism of fs log schemes, and 
assume that 
$f$ is a covering for the logarithmic flat topology. 
  Then, the underlying diagram of topological spaces of 
$$X\times_YX \rightrightarrows X \to Y$$
is exact, that is, $Y$ is the coequalizer of the left two arrows in the category of topological spaces.  
  (Here $X\times_YX$ is the fiber product in the category of fs log 
schemes.) 
\end{prop}

\begin{prop}\label{prop:kummer}
  Let $f\colon X \to Y$ be a morphism of fs log schemes. 

$(1)$  The following two conditions are equivalent. 

{\rm (i)} $f$ is of Kummer type. 

{\rm (ii)} \'Etale locally on $X$, 
there exists a covering $Y'\to Y$ for the log flat topology 
such that the log structure of $X\times_YY'$ is the inverse image 
of the log structure of $Y'$. 

$(2)$  Assume that 
$f$ is of Kummer type, $X$ is quasi-compact, and assume that 
we are given a chart $P \to M_Y$. 
  Then, globally on $X$, 
we can take as $Y'$ in {\rm (ii)} a log scheme over $Y$ of the 
type $Y\otimes_{\bZ[P]}\bZ[Q]$ endowed with the log structure 
associated to $Q \to \cO_{Y'}$, for some homomorphism of fs monoids 
$P \to Q$ of Kummer type.
\end{prop}

\noindent  {\emph Proof of $\ref{prop:open}$}. 
  Since classical fppf morphisms are open maps, it suffices to prove 
\ref{prop:open} in the following case: 
  $Y$ has a chart $P \to M_Y$, and there are an fs monoid $Q$ and a 
homomorphism $P \to Q$ of Kummer type such that 
$X=Y\otimes_{\bZ[P]}\bZ[Q]$ with the log structure associated to 
$Q\to \cO_{X}$. 
  Let $G$ be the group scheme $\Spec(\bZ[Q^{\gp}/P^{\gp}])$ over 
$\bZ$ with the trivial log structure. 
  Then, $G$ acts on $X$ over $Y$ by 
$$
\begin{matrix}
\cO_X = \cO_Y\otimes_{\bZ[P]}\bZ[Q]\to \cO_{G\times X} = 
\cO_Y \otimes_{\bZ[P]}\bZ[Q\oplus(Q^{\gp}/P^{\gp})]
\\ 
1\otimes a \mapsto 1 \otimes(a, a \bmod P^{\gp}) \qquad (a \in Q).
\end{matrix}
$$
  We have $G\times X \overset \cong \to X\times_YX; 
(g, x) \mapsto (x, gx)$. 
  (Recall that the fiber products are taken in the category of 
fs log schemes.) 

  Note that 
$f$ is a closed map (for the underlying morphisms of schemes 
of $f$ is finite) and surjective (indeed, since $\bZ[P]\to \bZ[Q]$ 
is injective and finite, $\Spec(\bZ[Q]) \to \Spec(\bZ[P])$ is 
surjective, and the surjectivity of a morphism of schemes is preserved 
by base changes). 
  Hence $Y$ has the quotient topology of the topology of $X$. 

  Let $U$ be an open set of $X$. 
  By the above remark, to see that $f(U)$ is open in $Y$, it is 
enough to show that $f^{-1}(f(U))$ is open in $X$. 
  But this last fact follows from \ref{lem:GxU} below and from the 
fact the action $G\times X \to X$ is fppf (since $G$ is fppf over 
$\bZ$). 

\begin{lem}\label{lem:GxU}
  With the notation as above, $f^{-1}(f(U))$ coincides with the image 
of $G\times U$ under the action $G\times X\to X$. 
\end{lem}

\begin{pf}
  This is reduced to the case where $Y$ is the spectrum of a field. 
\end{pf}

\medskip

\noindent 
{\emph Proof of $\ref{prop:exact}$.}
  By \ref{lem:4thpoint} (2), the sequence of the underlying sets is 
exact.  
  Hence, it is enough to show 
that the topology of $Y$ is the quotient 
topology of that of $X$. 
  Taking local charts, we may assume that 
there exists a chart $M_Y \leftarrow P \overset h \to 
Q \to M_X$ of $f$ with $h$ being of Kummer type 
such that the induced morphism of schemes $u\colon 
X \to Y\times_{\bZ[P]}\bZ[Q]$ is flat, surjective and locally 
of finite presentation. 
  Since 
$Y\times_{\bZ[P]}\bZ[Q] \to Y$ is closed and surjective, 
the topology of $Y$ is the quotient topology of that of 
$Y\times_{\bZ[P]}\bZ[Q]$. 
  Further, since $u$ is open and surjective, 
the topology of $Y\times_{\bZ[P]}\bZ[Q]$ is the quotient 
topology of that of $X$. 
  Thus we conclude that the topology of $Y$ is the quotient 
topology of that of $X$. 

\medskip

\noindent 
{\emph Proof of $\ref{prop:kummer}$.}
  Assume that the condition (ii) in (1) is satisfied.  
  Then,  $X\times_Y Y'\to X$ is of Kummer type by \ref{lem:4thpoint} (1) and  $X\times_YY' \to Y$ is of Kummer type, and hence $X \to Y$ is of Kummer type, that is, the condition (i) 
is satisfied.

  Next, assume that the condition (i) is satisfied. 
  Taking local charts for $f$, we see that, after localizing $X$, 
there is a positive integer $n$ such that the cokernel of 
$(f^*M_Y)^{\gp} \to M_X^{\gp}$ is killed by $n$.  
  We show that, under the assumption of the existence of such an $n$, 
globally on $X$, there exists a covering 
$Y' \to Y$ such that the log structure of $X \times_YY'$ is the 
inverse image of the log structure of $Y'$, which completes 
the proof of (1). 
  To show this, by \'etale localizing $Y$, we may and will 
assume that there is a chart $P \to M_Y$ of $Y$.
  We show that, under this further assumption of the existence 
of a chart of $Y$, there exists a covering of the type 
$Y\otimes_{\bZ[P]}\bZ[Q]$ as in (2). 
  This completes not only the proof of (1) but also the proof of (2) 
because there is an $n$ which kills 
the cokernel of $(f^*M_Y)^{\gp} \to M_X^{\gp}$ under the 
assumption of (2). 
  Let $h\colon P \to Q$ be a homomorphism $n\colon P \to P$ 
of Kummer type, let $Y'=Y\otimes_{\bZ[P]}\bZ[Q]$, and let 
$X'=X\times_YY'$. 
  It is enough to show that $f'{}^*M_{Y'} \to M_{X'}$ is an 
isomorphism, where $f'$ is the projection $X' \to Y'$. 
  Since $f'$ is of Kummer type, it is sufficient to show that 
$(f'{}^*M_{Y'})^{\gp} \to M_{X'}^{\gp}$ is an isomorphism, 
which is easily checked at stalks 
(cf.\ \cite{NC:loget} 2.1.1).

\section{Descent theory.\ I.}\label{sec:des1}

\bigskip 

  Concerning the descent theory for logarithmic flat topology, 
it seems a right philosophy is that {\emph morphisms descent,  
properties of morphisms also descent, 
but objects do not descent}. 
  We show 
in \S\ref{sec:des1} (resp.\ \S\ref{sec:des2}) that 
morphisms (resp.\ properties of morphisms) descent for the log 
flat topology (\ref{thm:des1}) (resp.\ (\ref{thm:des2})) and in \S\ref{sec:des3} 
that log objects descent for the classical flat topology (\ref{thm:des3}). 

\begin{thm}\label{thm:des1}
  Let $X$ be an fs log scheme, and let $Y$ be an fs log scheme over 
$X$. 
  Then, the functor 
$$\Mor_X(\ \ , Y)\colon T\mapsto \Mor_X(T, Y)$$
on $({\mathrm{fs}}/X)$ is a sheaf for the log flat topology.
\end{thm}

\begin{thm}
  Let $\Gmlog$ be the functor $T\mapsto \Gamma(T,M_T^{\gp})$ on 
$({\mathrm{fs}}/X)$. 
  Then $\Gmlog$ is a sheaf for the log flat topology.
\end{thm}

  We remark that in the circulated version (cf.\ Introduction), $\Gmlog$ was denoted by ${\mathbb G}_m^{\text{cpt}}$.

Since a sheaf for the log flat topology is also a sheaf for the 
log \'etale topology, the above results show that the functors 
$\Mor_X(\ \ ,Y)$ in \ref{thm:des1} and $\Gmlog$ are sheaves for 
the log \'etale topology.

\begin{para}\label{para:pf-des1}
  We prove Theorem \ref{thm:des1}.
  In this \ref{para:pf-des1}, we show that it is 
sufficient to prove that the functors 
\begin{equation}\tag{\thesubsection.1}\label{equ:OM}
T \mapsto \Gamma(T, \cO_T)\qquad {\text{and}} \qquad T \mapsto \Gamma(T, M_T)
\end{equation}
are sheaves for the log flat topology. 

  We may assume that $X=\Spec(\bZ)$ with the trivial log structure 
and 
$Y$ is affine as a scheme. 
  Further, we may assume that 
$Y$ has a chart $P \to M_Y$ (for this, see \cite{KKN4} \S 5). 

  Let $F, G, H$ be functors (fs$/\Spec(\bZ))\to($Set) defined by 
\begin{align*}
F(T)&= \{\text{ring homomorphisms }\Gamma(Y,\cO_Y) \to 
\Gamma(T,\cO_T)\}\\
G(T)&=\{\text{homomorphisms }P \to \Gamma(T, M_T)\} \\
H(T)&=\{\text{homomorphisms }P \to \Gamma(T, \cO_T)\}.
\end{align*}
  Then, $\Mor_X(\ \ ,Y)$ is the fiber product of $F\to H\leftarrow 
G$, where the first arrow is induced by $P\to\Gamma(Y, \cO_Y)$ and the 
second arrow is induced by $\Gamma(T, M_T) \to \Gamma(T, \cO_T)$. 
  Hence it is sufficient to prove that $F, G, H$ are sheaves. 

  Take a presentation 
\begin{align*}
&  \Gamma(Y, \cO_Y) = \bZ[T_i\ ; \ i \in I]/(f_j\ ; \ j \in J), \\
&  \bN^r \rightrightarrows \bN^s \to P\qquad ({\mathrm{exact}}). 
\end{align*}
  Then, $F(T)$ is the kernel of $\Gamma(T, \cO_T)^I \to \Gamma(T, \cO_T)^J$, 
and $G(T)$ and $H(T)$ are the equalizers of $\Gamma(T, M_T)^s 
\rightrightarrows \Gamma(T, M_T)^r$ and $\Gamma(T, \cO_T)^s 
\rightrightarrows \Gamma(T, \cO_T)^r$, respectively. 
  Hence it is sufficient to prove that the functors in \ref{equ:OM} 
are sheaves.
\end{para}

\begin{para}
  We prove that $T\mapsto\Gamma(T, \cO_T)$ 
is a sheaf for the log flat topology. 
  It is sufficient to prove 

\begin{sblem}\label{sblem:desO}
  Let $T$ be an fs log scheme which is affine as a scheme and which has a 
chart $P \to M_T$. 
  Let $Q$ be an fs monoid, let $P \to Q$ be a homomorphism of Kummer type, 
let $T' = T \otimes_{\bZ[P]}\bZ[Q]$ which we endow with the log structure associated 
to $Q \to \cO_{T'}$, and let $T'' = T' \times_TT'$ be the fiber product in the 
category of fs log schemes. 
  Then 
\begin{equation}\tag{$*$} 
\Gamma(T, \cO_T) \to 
\Gamma(T', \cO_{T'})\rightrightarrows 
\Gamma(T'', \cO_{T''})
\end{equation}
is exact.
\end{sblem}

\begin{pf}
  Let $A=\Gamma(T, \cO_T)$.
  Then, the diagram $(*)$ is isomorphic to 
\begin{equation}\tag{$**$} 
A \overset {\alpha} \to 
A \otimes_{\bZ[P]}\bZ[Q]\underset{\beta_2}{\overset {\beta_1}\rightrightarrows}
A \otimes_{\bZ[P]}\bZ[Q \oplus (Q^{\gp}/P^{\gp})],
\end{equation}
where $\beta_1$ (resp.\ $\beta_2$) is the homomorphism of $A$-algebras which 
sends $1\otimes a$ ($a \in Q$) to $1 \otimes (a, 1)$ (resp.\ 
$1 \otimes (a, a \bmod P^{\gp}$)). 
  Let $s\colon A \otimes_{\bZ[P]}\bZ[Q] \to A$ be the homomorphism of $A$-modules 
which sends $1\otimes a$ ($a \in Q$) to $a$ if $a \in P$, and to $0$ if $a \not\in P$. 
  Let $\iota\colon 
A \otimes_{\bZ[P]}\bZ[Q \oplus (Q^{\gp}/P^{\gp})] \to A \otimes_{\bZ[P]}\bZ[Q]$ 
be the  homomorphism of $A$-modules 
which sends $1\otimes (a, b)$ ($a \in Q, b \in Q^{\gp}/P^{\gp}$) to $1 \otimes a$ 
if $b \not=0$, and to $0$ if $b=0$. 
  (Here we use the fact that $P \to Q$ is of Kummer type to show that 
$s$ is 
well-defined.) 
  Then, $s\circ \alpha$ is the identity map of $A$, and $\alpha \circ s + \iota 
\circ (\beta_2 -\beta_1)$ is the identity map of $A\otimes_{\bZ[P]}\bZ[Q]$.
  This proves the exactness of $(**)$. 
\end{pf}
\end{para}

\begin{para}
  We prove that $T \mapsto \Gamma(T, M_T)$ is a sheaf for the log flat topology. 

  Let $T' \to T$ be a classical fppf covering (recall that this means 
that the log structure of $T'$ is the inverse image of that of $T$ and the 
underlying morphism of schemes of $T' \to T$ is an fppf covering), and 
let $T'' = T'\times_TT'$. 
  Then 
$$\Gamma(T, \cO_T^{\times}) \to 
\Gamma(T', \cO_{T'}^{\times})\rightrightarrows 
\Gamma(T'', \cO_{T''}^{\times})$$ 
is exact (a classical result). 
  Further, 
$$\Gamma(T, M_T/\cO_T^{\times}) \to 
\Gamma(T', M_{T'}/\cO_{T'}^{\times})\rightrightarrows 
\Gamma(T'', M_{T''}/\cO_{T''}^{\times})$$ 
is exact. 
  To see it, taking local charts, we may assume that $M_T/\cO_T^{\times}$ 
is the inverse image of a sheaf on the Zariski site, and use the fact that 
$M_{T'}/\cO_{T'}^{\times}$ and $M_{T''}/\cO_{T''}^{\times}$ are 
the inverse images of the sheaf $M_T/\cO_T^{\times}$ and the 
fact that the diagram of 
topological spaces 
$$T'' \rightrightarrows T' \to T$$ 
is exact. 
  We can deduce from these two exactness 
that $T \mapsto \Gamma(T, M_T)$ is a sheaf for the classical flat 
topology.

  Since $T \mapsto \Gamma(T, M_T)$ commutes with the limit with 
respect to the inverse system of the \'etale neighborhoods of a geometric 
point of the underlying scheme of $T$, 
it is sufficient to prove the following. 

\begin{sblem}
  Let $T, T', T'', P \to M_T$ and $P \to Q$ be as in the hypothesis of Lemma 
$\ref{sblem:desO}$, and assume $T=\Spec (A)$ for a local ring $A$, $P \overset 
\cong \to (M_T/\cO_T^{\times})_{\overline t}$, where $t$ is the 
closed point of 
$T$, and $Q$ has no torsion. 
  Then, 
$$\Gamma(T, M_T) \to \Gamma(T', M_{T'}) \rightrightarrows \Gamma(T'', M_{T''})$$ 
is exact.
\end{sblem}

\begin{pf}
  Let $A' = \Gamma(T', \cO_{T'}) = A \otimes_{\bZ[P]}\bZ[Q]$, 
$A'' = \Gamma(T'', \cO_{T''}) = A \otimes_{\bZ[P]}\bZ[Q \oplus (Q^{\gp}/P^{\gp})]$.
  By \ref{sblem:desO}, $A \to A' \rightrightarrows A''$ is exact.
  Let $I$ (resp.\ $I'$, resp.\ $I''$) 
be the ideal of $A$ (resp.\ $A'$, resp.\ $A''$) generated by the image of $P 
\setminus \{1\}$ (resp.\ $Q \setminus \{1\}$, resp.\ $Q \setminus \{1\}$), 
and let $V$ (resp.\ $V'$, resp.\ $V''$) 
be the subgroup of $A^{\times}$ (resp.\ $(A')^{\times}$, resp.\ $(A'')^{\times}$) 
consisting of elements which are congruent to $1$ modulo $I$ (resp.\ $I'$, resp.\ $I''$). 
  Since $A/I \overset \cong \to A'/I'$, we see that $V \to V' \rightrightarrows 
V''$ 
is exact. 
  It remains to show that 
\begin{equation}\tag{$*$} 
\Gamma(T, M_T)/V \to \Gamma(T', M_{T'})/V' \rightrightarrows \Gamma(T'', M_{T''})/V''
\end{equation}
is exact.
  Consider the exact diagram 
\begin{equation}\tag{$**$}
P \oplus (A/I)^{\times} \to Q \oplus (A/I)^{\times} 
\underset{\beta_2}{\overset {\beta_1}\rightrightarrows} Q
\oplus \{(A/I)[Q^{\gp}/P^{\gp}]\}^{\times},
\end{equation}
where $\beta_1$ (resp.\ $\beta_2$) is the map $(a,u) \mapsto (a,u)$ 
(resp.\ $(a, au)$).
  The exactness of this follows from the exactness of $P\to Q \rightrightarrows 
Q \oplus (Q^{\gp}/P^{\gp})$. 
  Now there is a natural homomorphism of diagrams from $(**)$ to $(*)$, 
which induces an isomorphism on the first terms, an isomorphism on 
the second terms, and an injection on the last terms.
  Since $(**)$ is exact, $(*)$ is also exact. 
\end{pf}
\end{para}

\begin{para}
  We prove that $\Gmlog$ is a sheaf for the log flat topology.
  Let $T' \to T$ be a log flat morphism of Kummer type which is surjective and locally of 
finite presentation as a morphism of schemes.
  Assume that $T$ has a chart $P \to M_T$ and is quasi-compact.
  Let $T''=T'\times_TT'$ (the fiber product in the category of fs log schemes).
  Then, $\Gamma(T, M_T^{\gp})$ (resp.\ $\Gamma(T', M_{T'}^{\gp})$, 
resp.\ $\Gamma(T'', M_{T''}^{\gp}))$ is isomorphic to $\underset a \varinjlim 
\Gamma(T, a^{-1}M_T)$ (resp.\ $\underset a \varinjlim 
\Gamma(T', a^{-1}M_{T'})$, resp.\ $\underset a \varinjlim 
\Gamma(T'', a^{-1}M_{T''})$), where $a$ ranges over all elements of $P$. 
  (This follows from the fact that $T'$ and $T''$ are of Kummer type over $T$.) 
  Hence the exactness of 
$\Gamma(T, M^{\gp}_T) \to \Gamma(T', M^{\gp}_{T'}) \rightrightarrows \Gamma(T'', M^{\gp}_{T''})$ 
is reduced to the exactness of 
$\Gamma(T, M_T) \to \Gamma(T', M_{T'}) \rightrightarrows \Gamma(T'', M_{T''})$.
  This implies that $\Gmlog$ is a sheaf for the log flat topology. 
\end{para}
 
\section{Cohomology.}\label{sec:coh}
  For an fs log scheme $X$, $X^{\cl}_{\fl}$ denotes the category $(\fs/X)$ endowed 
with the classical fppf topology.
  That is, $(U_i \to T)_i$ is a covering in $X^{\cl}_{\fl}$ means that $T$ is an 
object of $(\fs/X)$, the log structure of $U_i$ is the inverse image of $M_T$ for any 
$i$, the underlying scheme of $U_i$ is flat and locally of finite presentation over 
that of $T$, and the images of $U_i$ in $T$ cover $T$ set theoretically.

  The aim of this section is to prove 

\begin{thm}\label{thm:coh}
  Let $X$ be an fs log scheme and assume that $X$ is locally Noetherian as a scheme.
  Let $\ep\colon X^{\log}_{\fl} \to X^{\cl}_{\fl}$ be the canonical morphism of 
sites.
  Let $G$ be a commutative group scheme over the underlying scheme of $X$ 
satisfying either one of the following two conditions.

$(\mathrm i)$ $G$ is finite flat over the underlying scheme of $X$.
 
$(\mathrm{ii})$ $G$ is smooth and affine over the underlying scheme of $X$.

  We endow $G$ with the inverse image of the log structure of $X$. 
  Then we have a canonical isomorphism 

$$R^1\ep_*G \cong \underset {n \not=0} \varinjlim \cHom(\bZ/n (1), G) 
\otimes_{\bZ}(\Gmlog/\Gm),$$
where $n$ ranges over all non-zero integers and the inductive limit is 
taken with respect to the canonical projections 
$\bZ/mn(1) \to \bZ/n(1)$. 

  Here $\Gm$ is the functor $T \to \Gamma(T, \cO_T^{\times})$ on $(\fs/X)$ 
and $\bZ/n(1)=\Ker(n \colon \Gm \to \Gm)$ ($n\not=0$). 
  The quotient $\Gmlog/\Gm$ here is taken in the categories of sheaves 
on $X_{\fl}^{\cl}$.
\end{thm}

   The following Kummer sequence for $\Gmlog$ on $X_{\fl}^{\log}$ will 
be a starting point for the proof of \ref{thm:coh}. 

\begin{prop}
\label{prop:kummerseq}
  Let $X$ be an fs log scheme.
  Then the sequences of sheaves on $X_{\fl}^{\log}$ (resp.\ $X_{\et}^{\log}$) 
\begin{align*}
 0 & \to \bZ/n(1) \to \Gm \overset n \to \Gm \to 0 \\
 0 & \to \bZ/n(1) \to \Gmlog \overset n \to \Gmlog \to 0 
\end{align*}
are exact for any non-zero integer $n$ (resp.\ for any integer $n$ which 
is invertible on $X$). 
\end{prop}

\begin{pf}
  Similarly as in \cite{KN} 2.3. 
\end{pf}

\begin{para}
  Let $X$ be an fs log scheme and let $G$ be a sheaf of abelian groups on 
$X_{\fl}^{\log}$.
  We define a canonical homomorphism of sheaves on $X_{\fl}^{\log}$ 
\begin{equation*}\label{equ:R1}\tag{\thesubsection.1}
\varinjlim\cHom(\bZ/n(1), G)\otimes(\Gmlog/\Gm) \to R^1\ep_*G
\end{equation*}
as follows.
  Let $h$ be a local section of $\cHom(\bZ/n(1), G)$. 
  By the Kummer sequence 
$$
 0 \to \bZ/n(1) \to \Gmlog \overset n \to \Gmlog \to 0 
$$ 
on $X_{\fl}^{\log}$, we have 
$$\Gmlog = \ep_*\Gmlog \overset \delta \to R^1\ep_*(\bZ/n(1)) \overset h 
\to R^1\ep_*G,$$
where $\delta$ is the connecting homomorphism.
  The map $\delta$ kills $\Gm$, for $n \colon \Gm \to \Gm$ is surjective 
on $X_{\fl}^{\cl}$. 
  Thus we obtain the map \ref{equ:R1}.
\end{para}

\begin{para}
  We show that the case (i) of Theorem \ref{thm:coh} follows from the case 
(ii) of Theorem \ref{thm:coh}.
  Let $G$ be a finite flat commutative group over the underlying scheme 
of $X$. 
  Let $G^*$ be the Cartier dual of $G$, and let $L=\cMor(G^*, \Gm)$, where 
$\cMor$ means the sheaf of the morphisms of sheaves of sets.
  Then, we have an exact sequence 
$$0 \to G \to L \to L' \to 0, \qquad L' = L/G$$ 
and $L, L'$ are affine and smooth over the underlying scheme of $X$.
  By endowing $L$, $L'$ with the inverse images of the $M_X$, we have 
$$0 \to R^1\ep_*G \to R^1\ep_*L \to R^1\ep_*L' \qquad \mathrm{exact}$$
(for $L = \ep_*L \to L' =\ep_*L'$ is surjective on $X_{\fl}^{\cl}$).
  Hence the bijectivity of \ref{equ:R1} for $G$ is reduced to the 
bijectivities of \ref{equ:R1} for $L$ and $L'$.
\end{para}

\begin{para}\label{Xn}
  Assume that we are given a chart $P \to M_X$ such that $P$ is fs and torsion free 
(such chart exists classical \'etale locally on $X$).
  We study $R^1\ep_*$ by using \v{C}ech cohomology.

  For $n \ge1$, let 
$$X_n =X\otimes_{\bZ[P]}\bZ[P^{1/n}]$$
with the log structure associated to $\bZ[P^{1/n}] \to \cO_{X_n}$, 
and let $X_{n,i}$ ($i \ge0$) be the fiber product of $i+1$ copies of 
$X_n$ over $X$ in the category of fs log schemes. 
  Note that $X_n \to X$ is a covering in $X_{\fl}^{\log}$.
  For a sheaf $G$ of abelian groups on $X_{\fl}^{\log}$, we have 
a \v{C}ech complex $C_{G,n}$ 
$$C_{G,n} \colon \Gamma(X_{n,0}, G) \to \Gamma(X_{n,1}, G) \to 
\Gamma(X_{n,2}, G) \to \cdots.$$
\end{para}

\begin{lem}
\label{lem:Cech}
  Assume $X = \Spec(A)$ as a scheme with $A$ a strict local ring 
(i.e.\ a Henselian local ring with separably closed residue field), 
and assume that $G$ is represented by a smooth commutative group 
scheme over the underlying scheme of $X$ endowed with the inverse 
image of the log structure of $X$. 
  Then 
$$\underset n \varinjlim H^1(C_{G,n}) \overset \cong \to 
H^1(X_{\fl}^{\log}, G).$$
\end{lem}

\begin{pf} 
  By the general theory of \v{C}ech cohomology, 
$\underset n \varinjlim H^1(C_{G,n}) \to
H^1(X_{\fl}^{\log}, G)$ is injective and the cokernel is 
embedded in 
$\underset n \varinjlim H^1((X_n)_{\fl}^{\log}, G)$. 
  We show 
$\underset n \varinjlim H^1((X_n)_{\fl}^{\log}, G)=0$. 
  Let $a$ be an element of $H^1((X_n)_{\fl}^{\log}, G)$. 
  Then there exists a quasi-compact fs log scheme $T$ and a surjective 
log flat morphism of Kummer type $T \to X_n$ such that $a$ 
vanishes in $H^1(T_{\fl}^{\log}, G)$. 
  By \ref{prop:kummer} (2), for some $m \not=0$, $T\times_{X_n}
X_{mn}\to X_{mn}$ is with the inverse image log structure and flat in the 
classical sense.
  Hence $a$ vanishes classical flat locally on $X_{mn}$, so the 
image of $a$ in $H^1((X_{mn})_{\fl}^{\log}, G)$ comes from 
$H^1((X_{mn})^{\cl}_{\fl}, G)$. 
  Since $G$ is smooth, we have 
$H^1((X_{mn})^{\cl}_{\fl}, G)=
H^1((X_{mn})^{\cl}_{\et}, G)$ (\cite{Gr} 11.7). 
  Since $X_{mn}$ is the disjoint union of a finite number of Spec 
of strict local rings, we have $H^1((X_{mn})^{\cl}_{\et}, G)=0$.
\end{pf}

\begin{para}\label{para:brace}
  For $n \ge1$, let $H_n$ be the classical commutative group 
scheme $\Spec(\Z[P^{\gp}/(P^{\gp})^n])= \cHom(P^{\gp}, \bZ/n(1))$ over $\Spec(\bZ)$.
  Then $H_n$ acts on $X_n$ over $X$, and $H_n \times X_n \overset 
\cong \to X_n \times_XX_n$ (the fiber product is taken in the 
category of fs log schemes). 
  Hence we have 
$$X_{n, i} \cong \underbrace{H_n\times \cdots \times H_n}_{i\ 
\mathrm{times}}\times X_n.$$
\end{para}

\begin{para}\label{para:delta}
  Let $G$ be a sheaf of abelian groups on $X_{\fl}^{\log}$, and 
let $F$ be the sheaf of abelian groups on $X_{\fl}^{\log}$ defined 
by $F(T) = \Gamma(T \times_XX_n, G)$. 
  Then $H_n$ acts on $F$. 
  By \ref{para:brace}, $C_{G,n}$ is isomorphic to the complex 
\begin{align*}\tag{$*$}
  \Mor(\{1\}, F) \overset {\delta_0} \to \Mor(H_n, F) 
\overset {\delta_1} \to \Mor(H_n\times H_n, F) \to \cdots,
\end{align*} 
where $\Mor$ means the set of morphisms as sheaves of sets on 
$X_{\fl}^{\log}$, 
$$\delta_0(x)=(\sigma \mapsto \sigma x -x), \qquad 
\delta_1(x) = ((\sigma, \tau) \mapsto x(\sigma \tau) - x(\sigma) - 
\sigma(x(\tau))),\ldots.$$
\end{para}

\begin{para}
  Let $G$ and $F$ be as in \ref{para:delta}.
  Assume that $X$ is strictly local as a scheme and $P\overset 
\cong \to (M_X/\cO^{\times}_X)_x$, where $x$ is the closed point of $X$. 

\noindent 
Consider the trivial action of $H_n$ on $G$.
  Then, the natural homomorphism $G\to F$ preserves the action of $H_n$ 
and hence we get a homomorphism $H^1(H_n, G) \to H^1(H_n, F)$. 
  Note 
$$H^1(H_n, G) = \Hom(H_n, G) = \Hom(\bZ/n(1), G) \otimes_{\bZ}P^{\gp}.$$
  The composite map 
$$\Hom(\bZ/n(1), G) \otimes_{\bZ} P^{\gp} \cong H^1(H_n, G) \to H^1(H_n, F) 
\cong H^1(C_{G, n}) \to H^1(X_{\fl}^{\log}, G)$$ 
sends $h \otimes a$ to the image of $a$ under 
$$P^{\gp} \to H^0(X^{\log}_{\fl}, \Gmlog) \to H^1(X^{\log}_{\fl}, \bZ/n(1)) 
\overset h \to H^1(X_{\fl}^{\log}, G),$$ 
where the second arrow is the connecting homomorphism of 
$$
 0 \to \bZ/n(1) \to \Gmlog \overset n \to \Gmlog \to 0 \qquad \text{(cf.\ 
(\ref{prop:kummerseq})).}
$$ 
\end{para}

  Together with \ref{lem:Cech}, the next proposition \ref{prop:coh} implies \ref{thm:coh} 
under the assumption on $X$ in \ref{prop:coh}. 

\begin{prop}\label{prop:coh}
  Assume that as a scheme, $X=\Spec (A)$ for a Noetherian complete local ring $A$ 
with separably closed residue field.
  Let $G$ be a smooth commutative group scheme over the underlying scheme of $X$ 
and endow $G$ with the inverse image of the log structure of $X$. 
  Assume $P \overset \cong \to (M_X/\cO_X^{\times})_x$, where $x$ is the closed 
point of $X$.
  Then, 
$$\Hom(\bZ/n(1), G)\otimes_{\bZ} P^{\gp} \overset \cong \to 
H^1(C_{G,n}) \qquad \text{for any $n\ge1$}.$$
\end{prop}

\begin{para}\label{para:artin}
  We prove \ref{prop:coh} in the case where $A$ is an Artinian local ring.
  Let $I$ (resp.\ $J$) be the ideal of $A$ (resp.\ $\cO_{X_n}$) generated 
by the image of $P \setminus \{1\}$ (resp.\ $P^{1/n} \setminus \{1\})$. 
  Then $I$ (resp.\ $J$) is a nilpotent ideal.
  We define a descending filtrations $\{\fil^i G\}_{i\ge0}$ on the $H_n$-module 
$G$ and $\{\fil^i F\}_{i\ge0}$ on the $H_n$-module $F$ by 
\begin{align*}
(\fil^i G)(T)&=\Ker(G(T) \to G(T\times_X\Spec(\cO_X/I^i))) \\
(\fil^i F)(T)&=\Ker(F(T) \to G(T\times_X\Spec(\cO_{X_n}/J^i))).
\end{align*}
  By the nilpotence of $I$ (resp.\ $J$), we have $\fil^i =0$ for a sufficiently 
large $i$.
  Since $G$ is smooth, we have 

\begin{equation}\tag{\thesubsection.1}
\gr^i(G)(T)\cong \Lie(G)\otimes_A\Gamma(T, I^i\cO_T/I^{i+1}\cO_T),
\end{equation}
\begin{equation}\tag{\thesubsection.2}
\gr^i(F)(T)\cong \Lie(G)\otimes_A\Gamma(T, J^i\cO_T/J^{i+1}\cO_T)
\end{equation}
($\gr^i = \fil^i/\fil^{i+1}, i\ge1)$. 
  We have also 
\begin{equation}\tag{\thesubsection.3}\label{equ:gr0}
\gr^0(G)\overset \cong \to \gr^0(F).
\end{equation}
  By \ref{equ:gr0} and by the following \ref{lem:Cm}, we have 
\begin{multline*}
\Hom(\bZ/n(1), G) \otimes_{\bZ}P^{\gp} = H^1(H_n, G)  \\
\overset \cong \to 
H^1(H_n, \gr^0(G)) \overset \cong \to 
H^1(H_n, \gr^0(F)) \overset \cong \gets 
H^1(H_n, F)=H^1(C_{G,n}).
\end{multline*}
  This proves \ref{prop:coh} in the case $A$ is Artinian. 
\end{para}

\begin{lem}\label{lem:Cm}
  For any $i\ge1$ and any $m \ge1$, $H^m(H_n, \gr^i(G))$ and $H^m(H_n, \gr^i(F))$ 
are zero.
\end{lem}

\begin{pf}
  Since $\gr^i(G)$ and $\gr^i(F)$ are coherent, this follows from 
\cite{SGA3} I, 5.3.3.
\end{pf}

\begin{para}
  We prove Proposition \ref{prop:coh}.
  We denote the maximal ideal of $A$ by $m_A$.
  For $i \ge0$, define the abelian groups $D_i$ and $E_i$ by the exact sequences 
\begin{align*}
  0& \to G(A/m_A^i) \to G(X_{n,0} \otimes_AA/m_A^i) \to D_i \to  0\\
 0 & \to E_i \to G(X_{n,1} \otimes_AA/m_A^i) \to G(X_{n,2} \otimes_AA/m_A^i).
\end{align*}
  Then $D_i \subset E_i$ and $E_i/D_i$ is $H^1$ of the complex $C_{G,n}$ for 
$\Spec(A/m_A^i)$ which is endowed with the inverse image of the log structure of $X$. 
  Let $E = \underset i \varprojlim E_i$ and $D = \underset i \varprojlim D_i$.
  Since $G(A) \overset \cong \to \varprojlim G(A/m_A^i),$ 
$G(X_{n,k}) = \varprojlim G(X_{n,k}\otimes_AA/m_A^i),$ and since 
$G(A/m_A^{i+1}) \to G(A/m_A^{i})$ and $G(X_{n,0}\otimes_AA/m_A^{i+1}) \to 
G(X_{n,0}\otimes_AA/m_A^{i})$ are surjective by the smoothness of $G$, we have 
exact sequences 
\begin{align*}
  0& \to G(A) \to G(X_{n,0}) \to D \to 0 \\
 0 & \to E \to G(X_{n,1}) \to G(X_{n,2})
\end{align*}
and an isomorphism $E/D \cong \varprojlim E_i/D_i$.
Hence $E/D \cong H^1(C_{G,n})$.

  On the other hand, by \cite{SGA3} XV 1.6, 
$\Hom(\bZ/n(1), G)$ does not change when $X$ is replaced by $\Spec 
(A/m_A^i)$. 
  This and \ref{para:artin} show 
$$\Hom(\bZ/n(1), G)\otimes_{\bZ}P^{\gp}\overset \cong \to E_i/D_i\qquad \text{for any $i$}.$$
Hence 
$\Hom(\bZ/n(1), G)\otimes_{\bZ}P^{\gp}\overset \cong \to E/D \cong H^1(C_{G,n}).$
\end{para}

\begin{para}
  We prove \ref{thm:coh} in general. 
  We may assume $X=\Spec(A)$ as a scheme for a Noetherian strict local ring and let $G$ be a 
smooth affine commutative group scheme over $A$ endowed with the inverse image of $M_X$.

  Let $\hat X=\Spec(\hat A)$, where $\hat A$ is the completion of $A$, and endow $X$ with 
the inverse image of the log structure of $X$. 
  Since the problem is already solved for $\hat X$, and since $\Hom(\bZ/n(1), G)$ does not 
change when we replace $X$ by $\hat X$, it is sufficient to prove that $H^1(X_{\fl}^{\log}, 
G) \to H^1({\hat X}_{\fl}^{\log}, G)$ is injective.

  Let $\alpha$ be an element of $H^1(X_{\fl}^{\log}, G)$ which vanishes in 
$H^1({\hat X}_{\fl}^{\log}, G)$. 
  Since $\alpha$ vanishes classical fpqc locally, by the fpqc descent (here we use the fact that 
$G$ is affine), $\alpha$ is the class of a representable smooth affine $G$-torsor $Y$ over 
the underlying scheme of $X$ which is endowed with the inverse image of the log structure 
of $X$.
  Since $X$ is strictly local, $Y$ has an $X$-rational point and hence $Y$ is a trivial 
$G$-torsor.
  Hence $\alpha=0$.
\end{para}

\section{Hilbert 90.}\label{sec:Hilb}

\begin{thm}\label{thm:Hilb}
  Let $X$ be an fs log scheme whose underlying scheme is locally Noetherian.
  Then the canonical map from $H^1(X^{\cl}_{\et}, \Gmlog)$ to 
$H^1(X^{\log}_{\et}, \Gmlog)$ (resp.\ to $H^1(X^{\cl}_{\fl}, \Gmlog)$, 
resp.\ $H^1(X^{\log}_{\fl}, \Gmlog)$) is bijective.
\end{thm}

  Since $\Gmlog$ on the classical \'etale site is the direct image of $\Gmlog$ on 
the other sites, this theorem is equivalent to its local form 

\begin{cor}\label{cor:Hilb}
  Let $X$ be an fs log scheme whose underlying scheme is $\Spec$ of a Noetherian 
strict local ring.
  Then $H^1(X^{\cl}_{\fl}, \Gmlog)$, 
$H^1(X^{\log}_{\et}, \Gmlog)$, $H^1(X^{\log}_{\fl}, \Gmlog)$ are zero.
\end{cor}

  We prove \ref{cor:Hilb}.  
  Since the natural projections from $X^{\log}_{\fl}$ to 
$X^{\cl}_{\fl}$ and to $X^{\log}_{\et}$ send $\Gmlog$ to $\Gmlog$, 
it is enough to show $H^1(X^{\log}_{\fl}, \Gmlog)=0$.

  Consider the exact sequence $0 \to \Gm \to \Gmlog \to \Gmlog/\Gm \to 0$ 
on $X^{\log}_{\fl}$.  
  First, we show that the connecting map $H^0(X^{\log}_{\fl}, 
\Gmlog/\Gm) \to H^1(X^{\log}_{\fl}, \Gm)$ is surjective by using 
the result \ref{thm:coh} as follows. 
  Since $H^i(X^{\cl}_{\fl}, \Gm)=0$ for $i>0$ by \cite{Gr} 11.7, 
$H^1(X^{\log}_{\fl}, \Gm)=H^0(X^{\cl}_{\fl}, R^1\ep_*\Gm)$ and 
the latter is isomorphic to $H^0(X^{\cl}_{\fl}, (\bQ/\bZ) \otimes(\Gmlog/\Gm))$ 
by \ref{thm:coh}. 
  On the other hand, 
$H^0(X^{\log}_{\fl}, \Gmlog/\Gm)=H^0(X^{\cl}_{\fl}, \ep_*(\Gmlog/\Gm))
=H^0(X^{\cl}_{\fl}, \bQ \otimes(\Gmlog/\Gm))$.
  Hence, the cokernel of the above connecting map injects into 
$H^1(X^{\cl}_{\fl}, \Gmlog/\Gm)$, which is isomorphic to 
$H^1(X^{\cl}_{\et}, \Gmlog/\Gm)=0$ by \cite{Gr} 11.9. 

  Hence it is sufficient to prove 
$H^1(X^{\log}_{\fl}, \Gmlog/\Gm)=0$.
  To show this, we apply the same argument in the proof of \ref{lem:Cech}. 
  Take a chart by $P:=(M_{X}/\cO^{\times}_{X})_x$. 
  In the notation there, since $H^1((X_{mn})^{\cl}_{\fl}, \ep_*(\Gmlog/\Gm))
=H^1((X_{mn})^{\cl}_{\fl}, \bQ \otimes (\Gmlog/\Gm))
=H^1((X_{mn})^{\cl}_{\et}, \bQ \otimes (\Gmlog/\Gm))=0$ by \cite{Gr} 11.9 
again, 
it is enough to show 
$H^1(C_{\Gmlog/\Gm,n})=0$ for each fixed $n \ge1$.
  But the complex $C_{\Gmlog/\Gm,n}$ is easy to be described.  
  In fact, let $n'$ be the greatest divisor of $n$ which is invertible on $X$. 
  Let $R:=\Gamma(X_n, \Gmlog/\Gm)=
\bQ \otimes (P^{1/n})^{\gp}$ and $S:=\{1, \ldots, n'\}$. 
  Then, $X_{n,i}$ is a disjoint union of $(n')^i$ strict local fs log schemes, 
and the complex $C_{\Gmlog/\Gm,n}$ is isomorphic to 
$R \overset {\delta_0} \to \Map(S,R) \overset {\delta_1} 
\to \Map(S^2, R) \to \cdots$, 
where $(\delta_0(a))(i)=a$ for any $a \in R, i \in S$, and 
$(\delta_1(b))(i,j)=b(i)-b(j)$ for any $b \in \Map(S,R), i, j \in S$. 
  Hence $H^1(C_{\Gmlog/\Gm,n})=0$ as desired. 

\section{Locally free modules.}\label{sec:free}
\begin{para}
  In this section, for an fs log scheme $X$, we denote the sheaf $T \mapsto 
\Gamma(T, \cO_T)$ on $X_{\fl}^{\log}$ by $\cO_X$. 

  Consider the map 
\begin{equation*}\tag{\thesubsection.1}\label{equ:prod}
\overset n {\textstyle\prod} H^1(X_{\fl}^{\log}, \Gm) \to 
 H^1(X_{\fl}^{\log}, GL_n(\cO_X))
\end{equation*}
($\overset n \prod$ denotes the product of $n$ copies) induced by the 
embedding $\overset n \prod \Gm \overset \subset \to GL_n(\cO_X)$ as diagonal 
matrices. 
  If $X=\Spec(A)$ as a scheme for a Noetherian strict local ring, $H^1(X^{\log}_{\fl}, 
\Gm)$ is isomorphic to $(M^{\gp}_X/\cO_X^{\times})_x \otimes (\bQ/\bZ)$, where $x$ denotes 
the closed point of $X$ (\ref{thm:coh}).
\end{para}

\begin{thm}\label{thm:free}
  Let $X$ be an fs log scheme and assume $X = \Spec(A)$ as a scheme for a Noetherian 
strict local ring $A$.
  Then the map $\ref{equ:prod}$ induces an isomorphism of pointed sets 
$$H^1(X^{\log}_{\fl}, GL_n(\cO_X)) \cong \frak S_n \bs (\overset n {\textstyle\prod} 
((M_X^{\gp}/\cO_X^{\times})_x \otimes (\bQ/\bZ))).$$

  Here $\cO_X$ denotes the sheaf $T \to \Gamma(T, \cO_T)$ on $X_{\fl}^{\log}$, 
and $\frak S_n \bs$ means the quotient by the natural action of the symmetric 
group of degree $n$ on the product of $n$-copies.
\end{thm}

\begin{rem}
  For an fs log scheme $X$, 
$H^1(X_{\fl}^{\log}, GL_n(\cO_X))$ is identified with the set of isomorphism 
classes of $\cO_X$-modules on $X_{\fl}^{\log}$ which are locally free of 
finite rank.
  
  A way to obtain such modules is the following.
  Let $X$ be an fs log scheme having a chart $P \to M_X$, and take an fs monoid 
$Q$ and an integral homomorphism $P \to Q$ of Kummer type. 
  Let $Y$ be the scheme $X \otimes_{\bZ[P]}\bZ[Q]$ endowed with the log structure associated 
to $Q \to \cO_Y$. 
  Then if $f$ denotes the canonical morphism $Y_{\fl}^{\log} \to X_{\fl}^{\log}$, the 
$\cO_X$-module $f_*\cO_Y$ on $X_{\fl}^{\log}$ is locally free of finite rank.
  We have the following direct decomposition of $f_*\cO_Y$ into invertible 
modules.
  Let $H$ be the group scheme $\Spec (\bZ[Q^{\gp}/P^{\gp}])$ over $\bZ$.
  Then $H$ acts on $Y$ over $X$, and hence on $f_*\cO_Y$.
  For $a \in Q^{\gp}/P^{\gp}$, let $L_a$ be the part of $f_*\cO_Y$ on which 
$H$ acts via the character $H \to \Gm$ corresponding to $a$.
  Then 
$$f_*\cO_Y = \underset a \bigoplus L_a,$$ 
where $a$ ranges over all elements of $Q^{\gp}/P^{\gp}$ (this can be checked 
log flat locally).
  Let $a \in Q^{\gp}$ and let $m$ be a non-zero integer such that $a^m \in 
P^{\gp}$.
  Then the element of $H^1(X_{\fl}^{\log}, \Gm)$ corresponding to the 
invertible module $L_{a \bmod P^{\gp}}$ coincides with the image of $a^m$ 
under $H^0(X_{\fl}^{\log}, \Gmlog) \to H^1(X_{\fl}^{\log}, \bZ/m(1)) \to 
H^1(X_{\fl}^{\log}, \Gm)$, where the first arrow is the connecting map of the 
Kummer sequence $0 \to \bZ/m(1) \to \Gmlog \overset m \to \Gmlog \to 0$.
  If $X$ is as in \ref{thm:free} and $x$ is its closed point, this element of $H^1(X_{\fl}^{\log}, 
\Gm) = H^1(X_{\fl}^{\log}, GL_1(\cO_X))$ corresponds to $a^m \otimes m^{-1}$ of 
$(M^{\gp}_X/\cO_X^{\times})_x \otimes (\bQ/\bZ)$ in the isomorphism of \ref{thm:free}.
\end{rem}

\begin{cor}
  Let $X$ be as in $\ref{thm:free}$, and let $F$ be an $\cO_X$-module on $X_{\fl}^{\log}$ 
which is locally free of finite rank.
  Then $F$ is a direct sum of invertible $\cO_X$-modules on $X_{\fl}^{\log}$.
\end{cor}

\noindent 
  The method of proving of \ref{thm:free} is as follows. 
  Let $x$ be the closed point of $X$. 
  By the similar method to the proof of Theorem \ref{thm:coh}, we find 
$$H^1(X_{\fl}^{\log}, GL_n(\cO_X)) \cong \underset m \varinjlim\Hom(H_m, GL_n)/\sim,$$
where $H_m$ is the group scheme $\Hom((M_X^{\gp}/\cO_X^{\times})_x, \bZ/m(1))$ and $/\sim$ 
means the quotient set by the inner conjugation by elements of $GL_n(A)$.
  Theorem \ref{thm:free} can be deduced from this.
  
  This Theorem \ref{thm:free} is quoted in \cite{Ni} and a proof with details is given there. See 
   \cite{Ni} Theorem 3.22. 

  We prove some propositions concerning $\cO_X$-modules on log flat sites.
  Some results concerning $GL_n$-torsors and vector bundles on log flat sites are 
given in \S\ref{sec:torsor}.

\begin{propstar}\label{prop:zar}
  Let $X$ be an fs log scheme whose underlying scheme is affine, and let $F$ be an 
$\cO_X$-module satisfying the following condition$:$
  There is a covering $Y \to X$ in $X_{\fl}^{\log}$ such that the pullback of $F$ 
on $Y_{\fl}^{\log}$ is isomorphic to the module theoretic inverse image on 
$Y_{\fl}^{\log}$ of some quasi-coherent module on the small Zariski site of $Y$.
  Then, 
$$H^m(X_{\fl}^{\log}, F) = 0 \qquad \text{for any $m\ge1$}.$$
\end{propstar}

\noindent 
  The proof of \ref{prop:zar} is by the computation of \v Cech cohomology. 
  
 This Proposition \ref{prop:zar} is quoted in \cite{Ni} and a proof with details is given there. See \cite{Ni} Proposition 3.27. 

\begin{prop}\label{prop:cl}
  Let $X$ be an fs log scheme and let $0 \to F' \to F \to F'' \to 0$ be an 
exact sequence of $\cO_X$-modules on $X_{\fl}^{\log}$ which are locally free and 
of finite rank.

$(1)$ If the underlying scheme of $X$ is affine, this exact sequence splits.

$(2)$ $F$ is classical if and only if $F'$ and $F''$ are classical.

  Here we say an $\cO_X$-module $F$ on $X_{\fl}^{\log}$ which is locally free of 
finite rank is classical if the restriction $\overline F$ of $F$ to the small Zariski 
site of $X$ is locally free and $F$ is the module theoretic inverse image of $\overline 
F$ on $X_{\fl}^{\log}$.
\end{prop}

\noindent 
{\it Proof} of (1).
  Consider the exact sequence of cohomology groups associated to the exact sequence 
$$0 \to\cHom_{\cO_X}(F'', F') \to \cHom_{\cO_X}(F'', F) \to \cHom_{\cO_X}(F'', F'') \to 0$$
of sheaves on $X_{\fl}^{\log}$. 
  By \ref{prop:zar}, $H^1(X_{\fl}^{\log}, \cHom_{\cO_X}(F'', F'))=0$.
  Hence the identity map of $F''$ comes from $\Hom_{\cO_X}(F'', F)$, and hence we obtain 
the splitting.

\noindent 
{\it Proof} of (2).
  We may assume that the underlying scheme of $X$ is affine.
  Then, we are reduced to (1). 

\section{Descent theory.\ II.\ Logarithmic flat descent.}\label{sec:des2}
  In the circulated version (cf.\ Introduction), we have announced the following results (under 
some finiteness conditions). 
  The original proofs are so long and complicated, but 
for the first three cases of \ref{thm:des2} and for \ref{thm:ketdes}, Illusie--Nakayama--Tsuji (\cite{INT}) 
gave considerable short proofs based on a result of Olsson in \cite{O}. 
  So we do not include here the original proofs. 
  The last two cases of \ref{thm:des2} were also proved in Tani's \cite{T}. 

\begin{thmstar}\label{thm:des2}
  Let $f\colon X \to Y$ be a morphism of fs log schemes, and let $g\colon Y' \to Y$ be a 
covering in $Y_{\fl}^{\log}.$ 
  Let $X' = X\times_YY'$ (the fiber product in the category of fs log schemes) and let 
$f' \colon X' \to Y'$ be a morphism induced by $f$.

  Then 
$f$ is log \'etale (resp.\ log smooth, resp.\ log flat, 
resp.\ of Kummer type, resp.\ 
with finite underlying morphism of schemes) if and only if so is $f'$.
\end{thmstar}

\begin{thmstar}\label{thm:ketdes}
  Let $X' \overset g \to X \overset f \to Y$ be morphisms of fs log schemes, 
and assume that $g$ is surjective and of Kummer type.

  If $g$ and $f\circ g$ are log \'etale (resp.\ log smooth, resp.\ 
log flat), then $f$ is log \'etale (resp.\ log smooth, resp.\ log flat).
\end{thmstar}

\section{Descent theory.\ III.}\label{sec:des3}

  In the theory of fs log schemes, the descent theory for objects works for the 
classical fppf topology. 

\begin{thm}\label{thm:des3}
  Let $X$ be an fs log scheme and let $F$ be a sheaf on $X_{\fl}^{\cl}$.
  Assume that there is a covering $Y \to X$ in $X_{\fl}^{\cl}$ such that the 
inverse image of $F$ on $(\fs/Y)$ is represented by an fs log scheme $T$ over $Y$ 
which is of Kummer type over $Y$ and whose underlying scheme is affine over that of 
$Y$.
  
  Then, $F$ is represented by an fs log scheme $S$ over $X$ which is of Kummer type 
over $X$ and whose underlying scheme is affine over that of $X$. 
\end{thm}

  We remark that this theorem also follows from \cite{O} Appendix Corollary A.5. 

\begin{pf}
  Let $X, F, Y, T$ be as in the hypothesis of the theorem. 
  By the classical fppf descent, there exists a scheme $S$ over the scheme $X$ which 
is affine over the scheme $X$ and which is endowed with an isomorphism $T \cong 
S\times_XY$ of schemes over the scheme $Y$.
  Our task is to descend the log structure of $T$ to $S$.
  Let $T'$ be the fiber product $T \times_Y(Y\times_XY)$ in the category of fs log 
schemes, which represents the pullback of $F$ on $Y \times_XY$. 
  (As a scheme, $T'$ is the fiber product $T \times_ST$ in the category of schemes.) 
  Let $f\colon T \to S$ and $g \colon T' \to S$ be the canonical morphisms of 
schemes. 
  Note $f$ is faithfully flat and locally of finite presentation.
  Define the sheaf $M_S$ on $S_{\et}$ to be the equalizer of $f_*M_T \rightrightarrows 
g_*M_{T'}$ where the two arrows are induced by the first and the second projections 
$T' \to T$. 
  Then $M_S \to \cO_S$ is induced from the exactness of $\cO_S \to \cO_T 
\rightrightarrows \cO_{T'}$, and it is easy to see that $M_S$ is a log structure.
  Now \ref{thm:des3} is reduced to 
\begin{lem}\label{lem:fs}
  Let the notation be as above.
  Then, $M_S$ is an fs log structure and $M_T$ is the inverse image of $M_S$.
\end{lem}
  We deduce \ref{lem:fs} from 
\begin{lem}\label{lem:inverse}
  With the notation as above, $f^{-1}(M_S/\cO_S^{\times}) \overset \cong 
\to M_T/\cO_T^{\times}$. 
\end{lem}

\noindent 
{\it Proof} of \ref{lem:fs} assuming \ref{lem:inverse}. 
  Let $t \in T$ and $s = f(t) \in S$, and let $P=(M_S/\cO_S^{\times})_{\overline s}$. 
  Since $P$ is isomorphic to $(M_T/\cO_T^{\times})_{\overline t}$, it is an fs monoid 
(\ref{defn:fs}) (3).
  Take a homomorphism $h \colon P \to M_{S, \overline s}$ such that the composite 
$P \overset h \to M_{S, \overline s} \to P$ is the identity.
  Since $P$ is of finite presentation as a monoid, we can extend $h$ 
to a homomorphism $P \to M_S\vert_U$ for an \'etale neighbourhood $U$ of $\overline s$.
  Since $P \overset \cong \to (M_T/\cO_T^{\times})_{\overline t}$, 
$P \to M_T\vert_V$ is a chart for an \'etale neighbourhood $V$ of $\overline t$ in 
$f^{-1}(U)$.
  Let $W$ be the image of $V$ in $U$.
  Since $f$ is flat and locally of finite presentation, it is an open map 
(\cite{EGA4-2} Th\'eor\`eme 2.4.6), and hence $W$ is an open set of $U$. 
  Since the inverse image of $M_W/\cO_W^{\times}$ on $V$ coincides with $M_V/\cO_V^{\times}$, 
we have that $M_W$ coincides with the log structure associated to $P \to \cO_S$. 

  We prove \ref{lem:inverse} by dividing it into the two parts (1) (2) of 

\begin{lem}\label{lem:inverse2}
  Define the sheaf $N$ on $S_{\et}$ to be the equalizer of 
$f_*(M_T/\cO_T^{\times}) \rightrightarrows g_*(M_{T'}/\cO_{T'}^{\times})$. 
  Then{\rm :} 

$(1)$ $f^{-1}(N) \overset \cong \to M_T/\cO_T^{\times}$. 

$(2)$ $M_S/\cO_S^{\times} \overset \cong \to N$. 
\end{lem}

\noindent 
{\it Proof} of (1). 
  First note that all sheaves are for the classical \'etale topology here.
  Let $t \in T$, $a \in (M_T/\cO_T^{\times})_{\overline t}$ and let 
$s$ (resp.\ $x$) be the image of $t$ in $S$ (resp.\ $X$).
  We prove $a$ comes from $N_{\overline s}$.
  Since $T$ is of Kummer type over $X$, there exists $n \ge1$ such that 
$a^n$ comes from $b \in (M_X/\cO_X^{\times})_{\overline x}$. 
  Working classical \'etale locally on $S$ and on $X$, we may assume that $a$ comes 
from an element $\tilde a$ of $\Gamma(U, M_T/\cO_T^{\times})$ for 
some \'etale neighbourhood $U \to T$ of $\overline t\to T$ such that 
$U \to S$ is surjective, $b$ comes from an element $\tilde b$ 
of $\Gamma(X, M_X/\cO_X^{\times})$, and $\tilde a^n$ comes from $\tilde b$.
  (Here we used the fact $f$ is an open map.)

\begin{lem}\label{lem:power}
  For any $y \in T$, the image of $\tilde b$ in 
$(M_T/\cO_T^{\times})_{\overline y}$ is the $n$-th power of some element. 
\end{lem}

\begin{pf}
  Take $y' \in T$ which is in the image of $U$ in $T$ such that $f(y) = 
f(y')$, and take $z \in T'$ such that $p_1(z)=y$ and $p_2(z)=y'$ (here 
$p_i$ is the $i$-th projection $T' \to T$). 
  We have isomorphisms 
$$(M_T/\cO_T^{\times})_{\overline y}\overset \cong \to
(M_{T'}/\cO_{T'}^{\times})_{\overline z}\overset \cong \gets
(M_T/\cO_T^{\times})_{\overline {y'}}\quad \text{over\ }
(M_X/\cO_X^{\times})_{\overline x}.$$

  Now \ref{lem:power} follows from the fact that the image of $\tilde b$ 
in $(M_T/\cO_T^{\times})_{\overline {y'}}$ is 
an $n$-th power of the image of $\tilde a$.
\end{pf}

  By \ref{lem:power}, the image of $\tilde b$ in $M_T/\cO_T^{\times}$ 
is locally an $n$-th power.
  Since $M_T/\cO_T^{\times}$ is torsion free, the $n$-th root of $\tilde b$ 
in $M_T/\cO_T^{\times}$ exists globally on $T$, and it should coincide on 
$U$ with $\tilde a$.
  This $n$-th root of $\tilde b$ on $T$ is a global section of $N$ on $S$ by 
the uniqueness of $n$-th root in $M_{T'}/\cO_{T'}^{\times}$ of the image 
of $\tilde b$.
  This shows that $a$ comes from $N_{\overline s}$. This completes the 
  proof of (1) of \ref{lem:inverse2}. 

\medskip

\noindent 
{\it Proof} of \ref{lem:inverse2} (2).
  The problem is the surjectivity of $M_S \to N$.
  Let $h \in \Gamma(S, N)$.
  It is sufficient to prove that $h$ comes from $M_S$ \'etale 
locally on $S$.
  The inverse image of $h$ in $M_T$ under $M_T \to M_T/\cO_T^{\times}$ 
is an $\Gm$-torsor on $T$ which is endowed with descent data on $T'$.
  By the descent theory of line bundles, this $\Gm$-torsor descends to 
a $\Gm$-torsor $L$ on $S$.
  \'Etale locally on $S$, $L$ has a section which is regarded as a section 
of $M_S$ with image $h$ in $N$.
  Hence $M_S \to N$ is surjective.
\end{pf}

\section{Torsors.}\label{sec:torsor}
  The aim of this section is to prove

\begin{thm}\label{thm:torsor}
  Let $X$ be an fs log scheme which is locally Noetherian as a scheme.
  Let $G$ be a finite flat commutative group scheme over the underlying 
scheme of $X$, which we endow with the inverse image of $M_X$, and 
let $F$ be a $G$-principal homogeneous space in the category of sheaves 
on $X^{\log}_{\fl}$. 
  Then, $F$ is representable by an fs log scheme over $X$ which is log 
flat of Kummer type and whose underlying scheme is finite over that of 
$X$.
\end{thm}

\begin{para}
  Let $X$ and $G$ be as in the hypothesis of \ref{thm:torsor}.
  We denote by $H^1_{\mathrm r}(X^{\log}_{\fl}, G)$ the subset of $H^1(X^{\log}_{\fl}, 
G)$ consisting of elements whose corresponding $G$-principal homogeneous 
space on $X^{\log}_{\fl}$ is represented by an fs log scheme over $X$, 
which is log flat and of Kummer type and whose underlying scheme is finite 
over that of $X$.
  Theorem \ref{thm:torsor} states that 
$H^1_{\mathrm r}(X^{\log}_{\fl}, G)$ coincides with $H^1(X^{\log}_{\fl},G)$.  
  Note that 
$$H^1(X^{\cl}_{\fl}, G) \subset H^1_{\mathrm r}(X^{\log}_{\fl}, G) \qquad \text{in }
H^1(X^{\log}_{\fl}, G)$$ 
by the classical fppf descent theory.
\end{para}

\begin{lem}\label{lem:torsor}
  Let $X$ and $G$ be as in the hypothesis of $\ref{thm:torsor}$.

$(1)$ If $\chi_1 \in H^1_{\mathrm r}(X^{\log}_{\fl}, G)$ and 
$\chi_2 \in H^1(X^{\cl}_{\fl}, G)$, then 
$\chi_1+\chi_2 \in H^1_{\mathrm r}(X^{\log}_{\fl}, G)$.

$(2)$ Let $G'$ be a finite flat commutative group scheme over the 
underlying scheme of $X$, which we endow with the inverse image of the 
log structure of $X$, and let $h\colon G' \to G$ be an injective homomorphism.
  Then $h$ sends $H^1_{\mathrm r}(X^{\log}_{\fl}, G')$ into 
$H^1_{\mathrm r}(X^{\log}_{\fl}, G)$.
\end{lem}

\noindent 
{\it Proof} of (1).
  By \S\ref{sec:des3}, we may work classical fppf locally on $X$.
  Since $\chi_2$ vanishes classical fppf locally, we may assume $\chi_2=0$. 

\noindent 
{\it Proof} of (2).
  Let $F'$ be the $G'$-principal homogeneous space on $X^{\log}_{\fl}$ 
whose class belongs to $H^1_{\mathrm r}(X^{\log}_{\fl}, G')$. 
  Then the image of the class of $F'$ in $H^1(X^{\log}_{\fl}, G)$ is 
represented by the induced $G$-principal homogeneous space $F= 
G' \bs (G \times F')$, where $G'$ acts on $G \times F'$ by 
$(\sigma^{-1}, \sigma)$ $(\sigma \in G')$.
  We have a cartesian diagram 
$$
\begin{CD}
G \times F' @>>> F \\
@V{\mathrm pr}_1 VV @VVV \\
G @>>> G'\bs G.
\end{CD}
$$
  Since $G \to G' \bs G$ is fppf, the diagram implies that $F$ over $G' \bs G$ 
is represented classical fppf locally on $G' \bs G$, by a log flat fs log 
scheme of Kummer type whose underlying scheme is finite over the base.
  Hence by \S\ref{sec:des3}, $F \to G' \bs G$ is represented by an fs log scheme 
over $G' \bs G$ which is log flat and of Kummer type whose underlying scheme is 
finite over $G' \bs G$.\qed

\begin{lem}\label{lem:r}
  Let $X$ be an fs log scheme and let $n \ge1$.
  Then the image of the connecting map $H^0(X, M_X^{\gp}) \to 
H^1(X^{\log}_{\fl}, \bZ/n(1))$ of the Kummer sequence is contained in 
$H^1_{\mathrm r}(X^{\log}_{\fl}, \bZ/n(1))$.
\end{lem}

\begin{pf}
  Let $a \in H^0(X, M_X^{\gp})$, and let $F_a$ be the $\bZ/n(1)$-principal 
homogeneous space corresponding to the image of $a$ in 
$H^1(X^{\log}_{\fl}, \bZ/n(1))$.
  That is, $F_a(T) = \{b \in \Gamma(T, M_T^{\gp})\ ; \ b^n =a\}$ for any 
fs log scheme $T$ over $X$.
 By working classical \'etale locally on $X$, we may assume that there are a 
chart $P \to M_X$ and an element $\tilde a$ of $P$ whose image in 
$H^0(X, M_X^{\gp})$ is $a$.
  Let $L$ be the abelian group generated by $P^{\gp}$ and a letter $b$ which 
is subject to the relation $b^n = \tilde a$.
  Let $Q=\{x \in L\ ; \ x^n \in P\}$ and let $Y=X\otimes_{\bZ[P]}\bZ[Q]$ 
which is endowed with the log structure associated to $Q \to \cO_Y$.
  Then, $b \in F(Y)$ defines an isomorphism $\Mor_X(\ \ , Y) \overset \cong 
\to F_a$.
  Furthermore by the construction, $Y$ is log flat of Kummer type over $X$ 
and the underlying scheme of $Y$ is finite over that of $X$. 
\end{pf}

\begin{para}
  Now we prove \ref{thm:torsor}.
  By \S\ref{sec:des3} and by a limit argument, we may assume that as a scheme, 
$X=\Spec (A)$ for some Noetherian strict local ring $A$.
  Take a chart $P \to M_X$ which induces $P \overset \cong \to M_{X, x}/
\cO_{X,x}^{\times}$ where $x$ is the closed point of $x$.
  Let 
\begin{equation*}\tag{\thesubsection.1}\label{equ:torsor}
\underset n \varinjlim \Hom(\bZ/n(1), G)\otimes_{\bZ}P^{\gp} \to 
H^1(X^{\log}_{\fl}, G)
\end{equation*}
be the homomorphism which sends $h\otimes a$ ($h\in \Hom(\bZ/n(1), G)$, 
$a \in P^{\gp}$) to the image of $a$ under 
$P^{\gp} \to H^0(X, M_X^{\gp}) \to H^1(X^{\log}_{\fl}, \bZ/n(1)) \overset 
h \to H^1(X^{\log}_{\fl}, G)$.
  By \S\ref{sec:coh}, we have an isomorphism 
$$H^1(X^{\cl}_{\fl}, G) \oplus (\underset n \varinjlim \Hom(\bZ/n(1), G) 
\otimes_{\bZ}P^{\gp}) \overset \cong \to H^1(X^{\log}_{\fl}, G).$$
By \ref{lem:torsor} (1), it is enough to show that the image of 
$\underset n \varinjlim \Hom(\bZ/n(1), G) \otimes_{\bZ} P^{\gp} \to 
H^1(X^{\log}_{\fl}, G)$ is contained in $H^1_{\mathrm r}(X^{\log}_{\fl}, G)$.
  Let $G^{\mult}$ be the multiplicative part of $G$.
  Then, any homomorphism $\bZ/n(1) \to G$ factors through $G^{\mult}$.
  By \ref{lem:torsor} (2), we may assume $G$ is multiplicative.
  Then, $G \cong \underset i \oplus (\bZ/m_i\bZ)(1)$ for some finite family 
of non-zero integers $(m_i)_i$.
   Since a $G$-principal homogeneous space is the product of 
$\bZ/m_i(1)$-principal homogeneous spaces, we may assume $G=\bZ/m(1)$ for 
some $m\not=0$.
  Then, 
$$
\underset n \varinjlim \Hom(\bZ/n(1), G)\otimes_{\bZ}P^{\gp} \cong
P^{\gp}/(P^{\gp})^m,$$
and the map \ref{equ:torsor} is identified with the map 
$P^{\gp}/(P^{\gp})^m \to H^1(X^{\log}_{\fl}, \bZ/m(1))$ of the 
Kummer sequence.
  Hence we are reduced to \ref{lem:r}.
\end{para}

\begin{para}
  We add some remarks on $GL_n$-torsors and vector bundles on $X^{\log}_{\fl}$.
  We have seen in \S\ref{sec:free} that $\cO_X$-modules on $X^{\log}_{\fl}$ which 
are locally free of finite rank need not be classical ones.
  The following \ref{prop:GL} (resp.\ \ref{prop:vector}) says that representable 
$GL_n$-torsors (resp.\ vector bundles) on $X^{\log}_{\fl}$ {\it need not be} 
(resp.\ {\it need be}) classical ones.

  (Recall that the sheaf $\cO_X$ on $X^{\log}_{\fl}$ is defined by $T \to 
\Gamma(T, \cO_T)$.)
\end{para}

\begin{prop}\label{prop:GL}
  Let $F$ be a $GL_n(\cO_X)$-torsor in the category of sheaves on $X^{\log}_{\fl}$.
  Then $F$ is representable.
\end{prop}

\begin{prop}\label{prop:vector}
  Let $F$ be an $\cO_X$-module on $X^{\log}_{\fl}$ which is locally free of finite 
rank.
  Then $F$ is representable if and only if $F$ is a classical one.
\end{prop}

\noindent 
{\it Proof} of \ref{prop:GL}.
  By a limit argument we may assume that $X=\Spec(A)$ for a Noetherian strict local 
ring $A$.
  Let $\alpha \in H^1(X^{\log}_{\fl}, GL_n(\cO_X))$ and let $F$ be the 
$GL_n(\cO_X)$-torsor on $X^{\log}_{\fl}$ corresponding to $\alpha$.
  We show that $F$ is representable.
  Let $P\to M_X$ be a chart such that $P \overset \cong \to (M_X/\cO_X^{\times})_x$ 
where $x$ is the closed point of $X$.
  By \S\ref{sec:free}, $\alpha$ comes from an element $(\alpha_1, \ldots, \alpha_n)$ 
of $(P^{\gp} \otimes (\bQ/\bZ))^n$.
  Take $m_i \ge 1$ such that $m_i \alpha_i = 0$, and write $\alpha_i = a_i \otimes 
m_i^{-1}$, $a_i \in P^{\gp}$ and let $F_i$ be the $\bZ/m_i(1)$-torsor corresponding 
to the image in $H^1(X^{\log}_{\fl}, \bZ/m_i(1))$ of $a_i$ under the 
connecting map of $0 \to \bZ/m_i(1) \to \Gmlog \overset {m_i} \to \Gmlog \to 0$.
  Let $H=\underset{1 \le i \le n} \prod \bZ/m_i(1)$.
  Then the $GL_n(\cO_X)$-torsor $F$ is induced from the $H$-torsor 
$\underset{1 \le i \le n} \prod F_i$ by the diagonal embedding $H \to GL_n(\cO_X)$.
  Since the quotient $H \bs GL_n(\cO_X)$ is representable and $GL_n(\cO_X) \to 
H \bs GL_n(\cO_X)$ is fppf covering, the cartesian diagram 
$$
\begin{CD}
GL_n(\cO_X) \times \underset{1 \le i \le n} \prod F_i @>>> F \\
@VVV @VVV \\
GL_n(\cO_X) @>>> H \bs GL_n(\cO_X)
\end{CD}
$$
shows (\S\ref{sec:des3}) that $F$ is representable.

\noindent 
{\it Proof} of \ref{prop:vector}.
  The ``if" part is clear and so we consider the ``only if" part.
  We may assume $X =\Spec A$ as a scheme for a Noetherian strict 
local ring.
  Assume $F$ is representable but not classical.
  Then $F$ comes from an element $(a_i \otimes m_i^{-1})_{1 \le 
i \le n}$ of $((M_X^{\gp}/\cO_X^{\times})_x \otimes 
(\bQ/\bZ))^n$ ($n = \rank (F)$, $x$ is the closed point, $a_i 
\in (M_X^{\gp}/\cO_X^{\times})_x, m_i \in \bZ, m_i \ge 1$) 
such that $a_1$ is not an $m_1$-th power in 
$(M_X^{\gp}/\cO_X^{\times})_x$.

  Let $Y$ be the $GL_n(\cO_Y)$-torsor defined to be the subsheaf of 
$\overset n \prod F$ consisting of bases of $F$.
  Then $Y$ is representable (\ref{prop:GL}).
  We denote by the same letter $F$ (resp.\ $Y$) the fs log scheme 
over $X$ which represents $F$ (resp.\ $Y$).
  We prove the following $(*)$ and $(**)$.

  $(*)$  There is an open set $U$ of $\overset n \prod F$ which 
contains the image of the zero section $X \to \overset n \prod F$ 
such that $M_U$ coincides with the inverse image of $M_X$.

  $(**)$ For any $y \in Y$ lying over $x$, $(M_X^{\gp}/\cO_X^{\times})_{\overline x} 
\to (M_Y^{\gp}/\cO_Y^{\times})_{\overline y}$ is not bijective.

  These $(*)$ and $(**)$ imply $(U \cap Y) \times_X\{x\} = \varnothing$.
  But if $X' \to X$ is a covering in $X^{\log}_{\fl}$ such that the 
pullback of $F$ on $X'$ is classical and if $x'$ is a point of $X'$ 
lying over $x$, $(U \cap Y) \times_X\{x'\} \not= \varnothing$ because 
a classical vector bundle over a field is irreducible (so any two 
non-empty open sets intersect).

\noindent 
{\it Proof} of $(*)$.
  $\overset n \prod F$ is of Kummer type over $X$ since log flat locally 
on $X$, $F$ becomes a classical vector bundle with the inverse image of 
$M_X$.
  Let $y \in \overset n \prod F$ be a point in the image of the zero 
section $X \to \overset n \prod F$, and let $s$ be the image of $y$ in 
$X$.
  Then, the zero section defines $g\colon 
(M_Y^{\gp}/\cO_Y^{\times})_{\overline y} \to 
(M_X^{\gp}/\cO_X^{\times})_{\overline s}$ and the composite $g\circ h$ 
with the canonical map 
$h\colon (M_X^{\gp}/\cO_X^{\times})_{\overline s} \to 
(M_Y^{\gp}/\cO_Y^{\times})_{\overline y}$ is the identity map of 
$(M_X^{\gp}/\cO_X^{\times})_{\overline s}$.
  Since $h$ is of Kummer type, this means that $h$ is bijective.
  This proves $(*)$.

\noindent 
{\it Proof} of $(**)$.
  Since the pullback of $F$ on $Y$ is isomorphic to $\overset n \prod 
\cO_Y$, the image of $\alpha_i$ in 
$(M_Y^{\gp}/\cO_Y^{\times})_{\overline y}\otimes \bQ/\bZ$ vanishes 
for any $y\in Y$.
  Hence $a_1$ is an $m_1$-th power in $
(M_Y^{\gp}/\cO_Y^{\times})_{\overline y}$ for any $y \in Y$.

\section{Logarithmic fundamental groups.}\label{sec:pi1}
\begin{para}
  As is said in \S\ref{sec:des1}, the log flat descent for objects does 
not work well.
  An exception is the descent for log \'etale finite objects of Kummer type.

  For an fs log scheme $X$, let $\Et^{\log}(X)$ be the category of fs 
log schemes over $X$ which are log \'etale and of Kummer type over $X$ and 
whose underlying scheme is finite over that of $X$.
  Let $\lcf(X_{\et}^{\log})$ (resp.\ $\lcf(X_{\fl}^{\log})$) be 
the category of sheaves on $X_{\et}^{\log}$ (resp.\ $X_{\fl}^{\log}$) 
which are locally constant and finite.
\end{para}

\begin{thm}\label{thm:pi1}
  Let $X$ be an fs log scheme whose underlying scheme is locally 
Noetherian.

$(1)$ For an object $Y$ of $\Et^{\log}(X)$, the sheaf $\Mor_X(\ \ , Y)$ on 
$X^{\log}_{\et}$ (resp.\ $X_{\fl}^{\log}$) belongs to 
$\lcf(X_{\et}^{\log})$ (resp.\ $\lcf(X_{\fl}^{\log})$).

$(2)$ We have equivalences of categories 
$$\Et^{\log}(X) \overset{\simeq}\to \lcf(X^{\log}_{\et}) \overset{\simeq}\to 
\lcf(X_{\fl}^{\log}).$$

$(3)$ If $X$ is connected, there exists a profinite group $\pi_1^{\log}(X)$, 
which is unique up to isomorphism, such that the equivalent categories 
in $(2)$ are equivalent to the category of finite $\pi_1^{\log}(X)$-sets.

$(4)$ If $X=\Spec (A)$ as a scheme for a strict local ring $A$, 
we have an isomorphism 
$$\pi_1^{\log}(X) \cong \underset n \varprojlim \Hom((M_X^{\gp}/\cO_X^{\times})_x, 
\bZ/n(1))$$
where $n$ ranges over all integers which are invertible in $A$, $x$ is the closed point, and $\bZ/n(1)$ 
denotes $\Gamma(\Spec(A), \bZ/n(1))$.

  Here in $(3)$, for a profinite group $G$, a $G$-set means a discrete set endowed 
with a continuous action of $G$.
\end{thm}

  The log fundamental group is first treated in \cite{FK}. A large part of \ref{thm:pi1} is taken from \cite{FK}. 
  See also \cite{I}, \cite{S}, \cite{V}, \cite{Ho}.

\begin{rem}
  Let $X$ be a regular locally Noetherian scheme and let $D$ be a divisor on 
$X$ with normal crossings.
  Endow $X$ with the log structure defined by $D$ (that is, $M_X = \{f \in \cO_X\ 
; \ f$ is invertible outside $D\}$. 
  Then $\pi_1^{\log}(X)$ coincides with the tame fundamental group 
of Grothendieck-Murre $\pi_1(X, D)$ which controls finite \'etale schemes 
over $X\setminus D$ which are at worst tamely ramified at generic points of 
$D$.

\end{rem}

  In the rest of this section, we prove \ref{thm:pi1}. 

\begin{para}

  We prove \ref{thm:pi1} (1). 
  Note that $\Mor_X(\ \ , Y)$ is indeed a sheaf by \ref{thm:des1}. 
  Since a finite kummer log \'etale fs log scheme over an fs log scheme is kummer log \'etale locally a finite \'etale scheme in the classical sense endowed with the pullback 
log structure from the base, 
(1) is reduced to the nonlog corresponding result. 

\end{para}

\begin{prop}\label{pi1stl}

Let $X$ and $x$ be as in the assumption of  \ref{thm:pi1} (4) and let $G$ be a finite group. For an integer $n\geq 1$, let $H_n$ be the group scheme 
$\cH om((M^{\gp}_X/\cO_X^\times)_x, \Z/n(1))$.

$(1)$ Take a chart $P\to M_X$ such that $P\overset{\cong}\to (M_X/\cO_X^\times)_x$.
  For each integer $n\geq 1$, let $X_n\to X$ be as in \ref{Xn}. 
Let $T$ be a $G$-torsor on $X^{\log}_{\fl}$. Then there is an integer $n\geq 1$ which is invertible on $X$ such that the pullback of $T$ on $(X_n)^{\log}_{\fl}$ is a trivial $G$-torsor.

$(2)$ For Let $*=$ $\et$ or  $\fl$, $H^1(X^{\log}_*, G)$ is canonically isomorphic to  $\varinjlim_n \Hom(H_n(k), G)/\sim$,  where $n$ ranges over all integers $\geq 1$ which are invertible on $X$,  $k$ is the residue field of $x$, and $\sim$ is the $G$-conjugacy.

\end{prop}

\begin{pf} By 2.7 (2), if $T$ is a $G$-torsor on $X^{\log}_{\fl}$, there is an integer $n\geq 1$ such that the pullback of $T$ on
$(X_n)^{\log}_{\fl}$ is a $G$-torsor over $X_n$ in the classical sense with the inverse image of the log structure of $X_n$. Since this $G$-torsor over $X_n$ is finite \'etale over $X_n$ and $X_n$ is strict local, it is a trivial $G$-torsor. Thus we have proved a weaker version of (1) without the condition $n$ is invertible on $X$.

Hence $H^1(X^{\log}_{\fl}, G)$ is isomorphic to the inductive limit of the $H^1$- \v{C}ech cohomology  $\{g\in G(X_n\times_X X_n) \;|\; p_{13}^*(g) =p_{12}^*(g)p^*_{23}(g)\}/\sim$ of $(X_n/X,G)$, where $p_{ij}^*$ is the pullback by the $(i,j)$-th projection $p_{ij}: X_n \times_X X_n\times_X X_n \to X_n \times_X X_n$ and $\sim$ is the following equivalence relation.
 $g'\sim g$ if and only if $g'= p_1^*(u)gp_2^*(u)^{-1}$ for some $u\in G(X_n)$ where $p_i: X_n \times_X X_n\to X_n$ is the $i$-th projection. 
     
     For any Noetherian scheme $S$, $G(S)$ is identified with the set of all maps $\pi_0(S)\to G$ where $\pi_0(S)$ denotes the set of connected components of $S$. If $n'\geq 1$ denotes the largest divisor of $n$ which is invertible on $X$, we have $\pi_0(X_n \times_X X_n)\cong \pi_0(H_n \times X_n)= H_{n'}(k)$ and $\pi_0(X_n\times_X X_n \times_X X_n)\cong \pi_0(H_n \times H_n \times X_n)=H_{n'}(k)\times H_{n'}(k)$, and hence the $H^1$- \v{C}ech cohomology of $(X_n/X,G)$ is identified with $\Hom(H_{n'}(k), G)/\sim$ where $\sim$ is the conjugacy by $G$. Hence we obtain the isomorphism in (2) for $H^1(X^{\log}_{\fl}, G)$.
     This isomorphism is canonical (independent of the choice of the above chart $P\to M_X$) because the element of $H^1(X_{\fl}^{\log}, H_{n'}(k))$ corresponding to the identity map of $H_{n'}(k)$ comes from the canonical isomorphism in  \ref{thm:coh}.

This argument shows that the $H^1$-\v{C}ech cohomology of $(X_n/X, G)$ is the same as that of $(X_{n'}/X, G)$. This completes the proof of (1). 

 Since $H^1(X^{\log}_{\et}, G)\subset H^1(X^{\log}_{\fl}, G)$ and $X_n$ is log \'etale over $X$ if $n$ is invertible on $X$, we have $H^1(X^{\log}_{\et}, G)=H^1(X^{\log}_{\fl}, G)$. This completes the proof of (2). 
     \end{pf}

\begin{para}\label{pi12} Let $X$ be as in the assumption  of \ref{thm:pi1} (3). Let $*=$ \'et or fl and let ${\cal C}= \lcf(X_*^{\log})$. We prove that there exists a profinite group 
$\pi_{1,*}^{\log}(X)$ such that $\cal C$ is equivalent to the category of finite $\pi_{1,*}^{\log}(X)$-sets.

Let $k$ be an algebraically closed field with divisible integral log, let $n\geq 1$, and let $a:\Spec(k)\to X$ be a morphism (we call this morphism a base point). Let $\cal C$ be the category of sheaves on $X^{\log}_*$ which are locally constant and finite.  We prove that 
with the pullback functor by $a$, $\cal C$ is a Galois category in the sense of \cite{SGA1} Section V 5, which implies the above.

It is clear that $\cal C$ satisfies the condition (G1)--(G5) of Galois category in Section V 5 of \cite{SGA1}. We prove that $\cal C$ satisfies the remaining condition (G6) of the Galois category, that is, a morphism $S\to T$ of $\cal C$ is an isomorphism if the induced map  $S_a\to T_a$ is bijective. Let $I$ be the image of $S\to T$. Let $C$ be the complement of $I$ in $T$, so $T$ is the disjoint union of $I$ and $C$. Let $D$ be the complement of the diagonal image of $S$ in $S\times_I S$. Since $S_a\to T_a$ is bijective, $C_a$ and $D_a$ are empty. It is sufficient to prove that these $C$ and $D$ are empty. We have the vector bundles $\cO_C$ and $\cO_D$ on the site $X^{\log}_*$. Since $X$ is connected, $\cO_C$ and $\cO_D$ are of constant rank. Since they are $0$ at $a$, they are $0$. Hence $C$ and $D$ are empty.

\end{para}

\begin{para} 

We prove the equivalence $ \lcf(X^{\log}_{\et}) \overset{\simeq}\to 
\lcf(X_{\fl}^{\log})$ in \ref{thm:pi1} (2). We may assume that as a scheme, $X=\Spec(A)$ for a  strict local ring $A$ (and hence is connected). Then by \ref{pi12}, the first (resp. second) category is equivalent to the category of finite $\pi_{1,\et}^{\log}(X)$ (resp. $\pi^{\log}_{1,\fl}(X)$)-sets. By \ref{pi1stl}, we have 

(1) $\pi_{1,\et}^{\log}(X)= \pi_{1,\fl}^{\log}(X) \cong \varprojlim_n  \Hom((M_X^{\gp}/\cO_X^{\times})_x, 
\bZ/n(1))$ 

\noindent
where $n$ ranges over all integers which are invertible on $X$ and $x$ denotes the closed point of $X$. This proves  the part $ \lcf(X^{\log}_{\et}) \overset{\simeq}\to 
\lcf(X_{\fl}^{\log})$ of \ref{thm:pi1} (2).

Hence for any connected $X$, $\pi_{1,\et}^{\log}(X) = \pi_{1,\fl}^{\log}(X)$ and we denote this group by $\pi_1^{\log}(X)$.

We prove $\Et^{\log}(X) \overset{\simeq}\to \lcf(X^{\log}_{\et})$  in \ref{thm:pi1} (2).  We may assume that as a scheme, $X=\Spec(A)$ for a strict local ring $A$. By the above (1), it is sufficient to prove that for any commutative finite group $G$, any $G$-torsor over $X$ is represented by an object of $\Et^{\log}(X)$.  But this follows from Theorem \ref{thm:torsor}.

This completes the proof of \ref{thm:pi1}.

\end{para}

\bigskip

\noindent {\bf Correction to \cite{log}.}
  There are two corrections to \cite{log}. 
  The first is about Section 4 of \cite{log}, which was pointed by A.\ Ogus. 
  The definition of a morphism of weakly purely inseparable in \cite{log} (4.9) is strange and \cite{log} {\sc Proposition} (4.10) (2) and \cite{log} {\sc Lemma} (4.11) are false as they stand. 
  {\sc Proposition} (4.10) (1) is correct. We modify the paper as follows. 
  We delete \cite{log} (4.9), (4.10)(2), (4.11). 
  In \cite{log} {\sc Theorem} (4.12), we add the assumption that $f$ is of Cartier type, and $(X'', M'')=(X', M').$ 
 Note that in the rest of \cite{log}, {\sc Proposition} (4.12) is not used. In most applications, only morphisms of Cartier type are important.
 
  The next is about {\sc Proposition} (3.14) (4) of \cite{log}, which was pointed by Y.\ Nakkajima. 
  For this statement, it is necessary to assume that $f$ is integral. 
  Further, $I\cO_{\tilde X}$ there (which is not defined) should be replaced by $I\otimes_{\cO_Y} \cO_X$.

\noindent Kazuya Kato

\noindent 
Department of Mathematics
\\
University of Chicago
\\
5734 S.\ University Avenue
\\
Chicago, Illinois, 60637 \\
USA
\\
\noindent kkato@math.uchicago.edu

\bigskip 



\begin{thebibliography}{AMRT}




\bibitem[EGAIV]{EGA4-2}
  A.\ Grothendieck and J.\ Dieudonn\'e, \'Etude locale des sch\'emas et des morphismes de sch\'emas (EGA IV), Second partie, 
Publ.\ Math.\ Inst.\ Hautes \'Etude.\ Sci.\ 24 (1965), 5--231.



\bibitem[FK]{FK}
  K.\ Fujiwara and K.\ Kato, {Logarithmic \'etale topology theory}, preprint, 1995. 

\bibitem[G]{Gr}
  A.\ Grothendieck, Le groupe de Brauer III: exemples et compl\'ements, 
in Dix expos\'es sur la cohomologie des sch\'emas 
(A.\ Grothendieck and N.\ H.\ Kuiper, ed.), 
North-Holland, 1968.

\bibitem[Ha]{H}
  K.\ Hagihara, {Structure theorem of Kummer \'etale $K$-group}, 
J.\ of $K$-Theory 29 (2) (2003), 75--99.

\bibitem[Ho]{Ho}
  Y.\ Hoshi, {The exactness of the log homotopy sequence}, Hiroshima Math.\ J.\ 39 (1) (2009), 61--122.

\bibitem[I]{I}
  L.\ Illusie,
{An overview of the work of K.\ Fujiwara, K.\ Kato, and C.\ Nakayama on logarithmic \'etale cohomology}, Cohomologies $p$-adiques 
et applications Arithm\'etiques (II) (P.\ Berthelot, J.\ M.\ Fontaine, L.\ Illusie, K.\ Kato and M.\ Rapoport., \'ed.), 
Ast\'erisque 279 (2002), 271--322.

\bibitem[INT]{INT}
  L.\ Illusie, C.\ Nakayama, and T.\ Tsuji, 
On log flat descent, 
Proceedings of the Japan Academy, Series 89-A-1 (2013), 1--5.

\bibitem[$\mathrm K_1$]{log}
  K.\ Kato, {\rm Logarithmic structures of Fontaine-Illusie}, 
Algebraic analysis, geometry, and number theory (J.-I.\ Igusa, ed.), 
Johns Hopkins University Press, Baltimore, 1989, 
191--224.

\bibitem[$\mathrm K_2$]{K2}
  K.\ Kato, Logarithmic Dieudonn\'e theory, preprint.

\bibitem[KKN]{KKN4}
  K.\ Kato, T.\ Kajiwara, and C.\ Nakayama, {\rm Logarithmic abelian varieties, part IV: Proper models}, 
Nagoya Math.\ J.\ {\rm 219} (2015), 
9--63.

\bibitem[KN]{KN}
  K.\ Kato and C.\ Nakayama, {\rm Log Betti cohomology, log \'etale 
 cohomology, and log de Rham cohomology of log schemes over $\bold C$}, 
 Kodai Math.\ J. {\rm 22} (1999), 
161--186. 

\bibitem[KS]{KS}
  K.\ Kato and T.\ Saito, {On the conductor formula of Bloch}, 
Publ.\ Math., Inst.\ Hautes \'Etud.\ Sci.,  
100 (2004), 5--151.

\bibitem[M]{M}
  A.\ Moriwaki, Rigidity of morphisms for log schemes, preprint, 2005.

\bibitem[Na1]
{NC:loget}
  C.\ Nakayama, {Logarithmic \'etale cohomology}, Math.\ Ann.\ {308} 
(1997), 
365--404.

\bibitem[Na2]
{N:qs}
  C.\ Nakayama, {Quasi-sections in log geometry}, Osaka Math.\ J.\ {46} 
(2009), 
1163--1173.

\bibitem[Ni]
{Ni}
  W.\ Niziol, {$K$-theory of log schemes I}, Documenta Math.\ {13} 
(2008), 
505--551.

\bibitem[Og]{Og}
  A.\ Ogus, {Lectures on Logarithmic Algebraic Geometry}, 
Cambridge Studies in Advanced Mathematics 178, Cambridge University Press, 2018.

\bibitem[Ol]{O}
  M.\ C.\ Olsson, Logarithmic geometry and algebraic stacks, 
Ann.\ Scient.\ \'Ec.\ Norm.\ Sup., $4^{\rm e}$ s\'erie, 
36 (2003), 747--791.

\bibitem[S]{S}
  J.\ Stix, Projective anabelian curves in positive characteristic and descent theory for log-\'etale covers,  Dissertation (Rheinische Friedrich-Wilhelms-Universit\"at Bonn, 2002), Bonner Math.\ Schriften, 354, Universit\"at Bonn, Mathematisches Institut, 2002.
  
  
\bibitem[SGA1]{SGA1} 
 A. Grothendieck,
 {Rev\^etements \'etales et groupe fondamental}, Lecture Notes in Math. 224 (1971).



\bibitem[SGA3]
{SGA3}
  M.\ Demazure and A.\ Grothendieck with M.\ Artin, J.\ E.\ Bertin, 
P.\ Gabriel, M.\ Raynaud, J.\ P.\ Serre, 
{Sch\'emas en Groupes}, Lecture Notes in Mathematics {151--153} 
(1970).

\bibitem[T]{T}
  A.\ Tani, 
Log flat descent of finiteness, thesis, Tokyo Institute of Technology, 2014 
(in Japanese).

\bibitem[V]{V}
  I.\ Vidal,  Contributions \`a la cohomologie \'etale des sch\'emas et des log-sch\'emas,  Th\'ese, U. Paris-Sud, 2001.
\end{thebibliography}
\end{document}